\input amstex
\input amsppt.sty \magnification=1200
\NoBlackBoxes
\hsize=15.6truecm \vsize=22.2truecm
\def\q{\quad}
\def\qq{\qquad}
\def\qtq#1{\q\t{#1}\q}

\def\({\left(}
\def\){\right)}
\def\[{\left[}
\def\]{\right]}
\def\mod#1{\ (\text{\rm mod}\ #1)}
\def\t{\text}
\def\f{\frac}

\def\e{\equiv}
\def\a{\alpha}
\def\b{\binom}

\def\ap{\langle a\rangle_p}

\def\sls#1#2{(\f{#1}{#2})}
\def\ag#1{\langle #1\rangle_p}
\def\ls#1#2{\big(\f{#1}{#2}\big)}
\def\Ls#1#2{\Big(\f{#1}{#2}\Big)}
\let \pro=\proclaim
\let \endpro=\endproclaim

\topmatter
\title {Jacobsthal sums, Legendre polynomials and binary quadratic forms}\endtitle
\author Zhi-Hong Sun \endauthor
\affil School of Mathematical Sciences, Huaiyin Normal University,
\\ Huaian, Jiangsu 223001, PR China
\\ Email: zhihongsun$\@$yahoo.com
\\ Homepage: http://www.hytc.edu.cn/xsjl/szh
\endaffil

 \nologo \NoRunningHeads
 \nologo \NoRunningHeads

\abstract{Let $p>3$ be a prime and $m,n\in\Bbb Z$ with $p\nmid mn$.
Built on the work of Morton, in the paper we prove the uniform
congruence:
$$\aligned &\sum_{x=0}^{p-1}\Big(\frac{x^3+mx+n}p\Big)
\\&\e\cases  -(-3m)^{\frac{p-1}4}
\sum_{k=0}^{p-1}\binom{-\frac 1{12}}k\binom{-\frac 5{12}}k
(\frac{4m^3+27n^2}{4m^3})^k\pmod p&\t{if $4\mid p-1$,}
\\\frac{2m}{9n}(\frac{-3m}p)(-3m)^{\frac{p+1}4}
\sum_{k=0}^{p-1}\binom{-\frac 1{12}}k\binom{-\frac 5{12}}k
(\frac{4m^3+27n^2}{4m^3})^k\pmod p&\text{if $4\mid p-3$,}
\endcases\endaligned$$
 where $(\frac ap)$ is the Legendre symbol.
 We also establish  many congruences for $x\pmod p$, where
 $x$ is given by $p=x^2+dy^2$ or $4p=x^2+dy^2$, and
 pose some conjectures on supercongruences modulo $p^2$
concerning binary quadratic forms.
\par\q
}
 \endabstract
  \footnote"" {The author was
supported by the Natural Sciences Foundation of China (grant no.
10971078).}
\endtopmatter
\document
\subheading{1. Introduction}

\par Let $\{P_n(x)\}$ be the Legendre polynomials given by
$$P_0(x)=1,\ P_1(x)=x,\ (n+1)P_{n+1}(x)=(2n+1)xP_n(x)-nP_{n-1}(x)\ (n\ge
1).$$  It is well known that (see [G, (3.132)-(3.133)])
$$P_n(x)=\f
1{2^n}\sum_{k=0}^{[n/2]}\b nk(-1)^k\b{2n-2k}nx^{n-2k} =\f 1{2^n\cdot
n!}\cdot\f{d^n}{dx^n}(x^2-1)^n,\tag 1.1$$ where $[a]$ is the
greatest integer not exceeding $a$.
 From (1.1)
we see that
$$P_n(-x)=(-1)^nP_n(x), \q P_{2m+1}(0)=0\qtq{and}P_{2m}(0)
=\f{(-1)^m}{2^{2m}}\b{2m}m.\tag 1.2$$
 We also have the following
formula due to Murphy ([G, (3.135)]):
$$P_n(x)=\sum_{k=0}^n\b nk\b{n+k}k\Ls{x-1}2^k=
\sum_{k=0}^n\b{2k}k\b{n+k}{2k}\Ls{x-1}2^k.\tag 1.3$$ We remark that
$\b nk\b{n+k}k=\b{2k}k\b{n+k}{2k}$.
\par Let $\Bbb Z$ be the set of
integers. For a prime $p$ let $\Bbb Z_p$ be the set of rational
numbers whose denominator is coprime to $p$. Based on the work of
Morton ([BM],[M]), in [S3-S6] the author proved that for a prime
$p>3$ and $t\in\Bbb Z_p$,
$$\align
&P_{\f{p-1}2}(t)\e
-\Ls{-6}p\sum_{x=0}^{p-1}\Ls{x^3-3(t^2+3)x+2t(t^2-9)}p\mod p,\tag
1.4 \\&P_{[\f p4]}(t) \e
-\Ls{6}p\sum_{x=0}^{p-1}\Ls{x^3-\f{3(3t+5)}2x+9t+7}p\mod
 p,\tag 1.5
 \\&P_{[\f p3]}(t)\e -\Ls
 p3\sum_{x=0}^{p-1}\Ls{x^3+3(4t-5)x+2(2t^2-14t+11)}p\mod p,\tag 1.6
\\&P_{[\f p6]}(t)
\e -\Ls 3p\sum_{x=0}^{p-1}\Ls{x^3-3x+2t}p\mod p,\tag 1.7\endalign
$$ where $\ls ap$ is the Legendre symbol.
\par For a prime $p>3$ let $F_p$ be the field of $p$ elements. For $m,n\in F_p$ we call the
curve $E_p:\ y^2=x^3+mx+n$ an elliptic curve provided that the
discriminant $4m^3+27n^2\not=0$ in $F_p$. Let $\#E_p(y^2=x^3+mx+n)$
be the number of points (including infinity) on the curve $E_p:\
y^2=x^3+mx+n$ over $F_p$. It is easily seen that (see for example
[S1])
$$\#E_p(y^2=x^3+mx+n)=p+1+\sum_{k=0}^{p-1}\Ls{x^3+mx+n}p.\tag 1.8$$
The famous inequality due to Hasse states that
$$\Big|\sum_{k=0}^{p-1}\Ls{x^3+mx+n}p\Big|\le 2\sqrt p.$$
It is important and difficult to evaluate the general Jacobsthal sum
$\sum_{k=0}^{p-1}\sls{x^3+mx+n}p$. So far, the sum is only
determined by using Deuring's theorem when the corresponding
elliptic curve has complex multiplication (see [LM],[I]). In 2006,
Morton[M] established a congruence for
$\sum_{k=0}^{p-1}\sls{x^3+mx+n}p\mod p$ by using Jacobi polynomials.
\par  Let $p>3$ be a prime and $m,n\in\Bbb Z_p$ with $mn\not\e
0\mod p$. Built on the work of Morton, in the paper we prove the
uniform congruence:
$$\aligned &\sum_{x=0}^{p-1}\Ls{x^3+mx+n}p
\\&\e\cases  -(-3m)^{\f{p-1}4}
\sum_{k=0}^{p-1}\b{-\f 1{12}}k\b{-\f 5{12}}k
(\f{4m^3+27n^2}{4m^3})^k\mod p&\t{if $4\mid p-1$,}
\\\f{2m}{9n}\sls{-3m}p(-3m)^{\f{p+1}4}
\sum_{k=0}^{p-1}\b{-\f 1{12}}k\b{-\f 5{12}}k
(\f{4m^3+27n^2}{4m^3})^k\mod p&\t{if $4\mid p-3$.}
\endcases
\\&\e\cases  (-1)^{\f{p+1}2}\sls n2^{\f{p-1}6}
\sum_{k=0}^{p-1}\b{-\f 1{12}}k\b{-\f 7{12}}k
\sls{4m^3+27n^2}{27n^2}^k\mod p&\t{if $3\mid p-1$,}
\\ (-1)^{\f{p+1}2}\f 3m\sls 2n^{\f{p-5}6}
\sum_{k=0}^{p-1}\b{-\f 1{12}}k\b{-\f 7{12}}k
\sls{4m^3+27n^2}{27n^2}^k\mod p&\t{if $3\mid p-2$.}\endcases
\endaligned\tag 1.9$$
We note that the right sums in (1.9) have the good form
$\sum_{k=0}^{p-1}f(k)$.
\par \par For positive integers $a,b$ and $n$,
 if $n=ax^2+by^2$ for some $x,y\in\Bbb Z$, we briefly say that
$n=ax^2+by^2$. Let $p$ be a prime of the form $4k+1$ and so
$p=x^2+y^2$ with $x\e 1\mod 4$. In 1825, Gauss found the congruence
$$2x\e \b{\f{p-1}2}{\f{p-1}4}\mod p.$$
Similar congruences for $x\mod p$ with
$p=x^2+dy^2(d\in\{2,3,5,7,11\})$ were found by Jacobi, Eisenstein
and Cauchy, see [BEW]. For $d=3,7,11,19,43,67,163$ we know that for
any prime $p$ with $\sls pd=1$, there are unique positive integers
$x$ and $y$ such that $4p=x^2+dy^2$. In [R1,R2], [RPR], [PR], [PV]
and [JM], the $x$ was given by an appropriate cubic Jacobsthal sum.
For example,
$$
\sum_{x=0}^{p-1}\Ls{x^3-96\cdot 11x+112\cdot 11^2}p=\cases \sls
2p\sls u{11}u&\t{if $\sls p{11}=1$ and $4p=u^2+11v^2$,}\\0&\t{if
$\sls p{11}=-1$.}\endcases$$
\par In the paper, based on the work of Ishii[I], by considering
the sums
$$\sum_{k=0}^{[p/6]}\b{[\f p3]}k\b{[\f p6]}km^k,\
\sum_{k=0}^{[p/8]}\b{[\f p8]}k\b{[\f {3p}8]}km^k, \
\sum_{k=0}^{p-1}\b{-\f 13}k^2m^k, \ \sum_{k=0}^{p-1}\b{-\f
14}k^2m^k$$ we establish many congruences for $x\mod p$, where $p$
is an odd prime and $x\in\Bbb Z$ is given by $p=x^2+dy^2$ or
$4p=x^2+dy^2$. Here are four typical results:
$$\aligned&2x\e\sum_{k=0}^{(p-1)/6}\b{\f {p-1}3}k\b{\f {p-1}6}k(-4)^k\mod
p\qtq{for} p=x^2+15y^2,
\\&x\e \sum_{k=0}^{(p-1)/6}\b{\f {p-1}3}k\b{\f {p-1}6}k\f 1{(-16)^k}\mod
p\qtq{for} 4p=x^2+51y^2,
\\&2x\e \sum_{k=0}^{(p-1)/8}\b{\f {p-1}8}k\b{\f {3(p-1)}8}k\f 1{(-882^2)^k}\mod
p\qtq{for} p=x^2+37y^2\e 1\mod 8,
\\&2x\e \sum_{k=0}^{[p/8]}\b{[\f p8]}k\b{[\f {3p}8]}k\f 1{99^{4k}}\mod
p\qtq{for} p=x^2+58y^2.
\endaligned$$
We also pose many conjectures on congruences modulo $p^2$. For
example, we conjecture that $$y\e \f{910}{9801}
\sum_{k=0}^{[p/8]}\b{[\f p{8}]}k\b{[\f {3p}8]}k\f 1{99^{4k}} \mod p
\qtq{for}p=2x^2+29y^2.$$

 \subheading{2. General congruences involving
$\b ak\b{-\f 12-a}k$, $\b ak\b{a-\f12}k$ and $\b ak^2$}
\par  For a prime $p$ and $a\in\Bbb Z_p$ let $\ap\in\{0,1,\ldots,p-1\}$
be given by $a\e \ap\mod p.$ Let $P_n^{(\a,\beta)}(x)$ be the Jacobi
polynomial defined by
$$P_n^{(\a,\beta)}(x)=\f 1{2^n}\sum_{k=0}^n\b{n+\a}k\b{n+\beta}{n-k}
(x+1)^k(x-1)^{n-k}.\tag 2.1$$ It is known that (see [AAR, p.315])
$$P_{2n}(x)=P_n^{(0,-\f 12)}(2x^2-1)\qtq{and}
P_{2n+1}(x)=xP_n^{(0,\f 12)}(2x^2-1).\tag 2.2$$ From [B, p.170] we
know that
$$\align P_n^{(\a,\beta)}(x)&=\b{n+\a}n\sum_{k=0}^{\infty}
\f{(-n)_k(n+\a+\beta+1)_k}{(\a+1)_k\cdot k!}\Ls{1-x}2^k
\\&=\b{n+\a}n\sum_{k=0}^n\f{\b
nk\b{-n-\a-\beta-1}k}{\b{-1-\a}k}\Ls{x-1}2^k.\endalign$$ Thus,
$$P_n^{(0,\beta)}(x)=\sum_{k=0}^n\b nk\b{-n-\beta-1}k\Ls{1-x}2^k.\tag 2.3$$

\pro{Lemma 2.1} Let $n$ be a nonnegative integer. Then
$$\align &P_n(x)=x^n\sum_{k=0}^{[n/2]}\b n{2k}\b{-\f 12}k
\Big(\f 1{x^2}-1\Big)^k,
\\&P_{2n}(x)=\sum_{k=0}^n\b nk\b{-\f 12-n}k(1-x^2)^k,
\\&P_{2n+1}(x)=x\sum_{k=0}^n\b nk\b{-\f 32-n}k(1-x^2)^k.\endalign$$
\endpro
Proof. Since $\b{-\f 12}k=\b{2k}k\f 1{(-4)^k}$, by [G, (3.137)] we
have
$$P_n(x)=x^n\sum_{k=0}^{[n/2]}\b n{2k}\b{2k}k\Ls{x^2-1}{4x^2}^k
=x^n\sum_{k=0}^{[n/2]}\b n{2k}\b{-\f 12}k\Big(\f 1{x^2}-1\Big)^k.$$
From (2.2) and (2.3) we see that
$$\align P_{2n}(x)&=P_n^{(0,-\f 12)}(2x^2-1)
=\sum_{k=0}^n\b nk\b{-n+\f 12-1}k\Ls{1-(2x^2-1)}2^k
\\&=\sum_{k=0}^n\b nk\b{-\f 12-n}k(1-x^2)^k\endalign$$
and
$$\align P_{2n+1}(x)&=xP_n^{(0,\f 12)}(2x^2-1)
=x\sum_{k=0}^n\b nk\b{-n-\f 12-1}k\Ls{1-(2x^2-1)}2^k
\\&=x\sum_{k=0}^n\b nk\b{-\f 32-n}k(1-x^2)^k.\endalign$$
This proves the lemma.

\pro{Theorem 2.1} Let $p$ be an odd prime and let $m$ be a positive
integer such that $p\nmid m$. Then
$$\aligned&\sum_{k=0}^{[\f p{2m}]}\b{[\f p{2m}]}k\b{[\f{(m-1)p}{2m}]}k(1-t)^k
\\&\e t^{[\f p{2m}]}\sum_{k=0}^{[\f p{2m}]}\b{[\f p{m}]}{2k}\b{\f{p-1}2}k
\Big(\f 1t-1\Big)^k
\\&\e\cases P_{[\f pm]}(\sqrt t)\mod p&\t{if $2\mid [\f pm]$,}
\\P_{[\f pm]}(\sqrt t)/\sqrt t\mod p&\t{if $2\nmid [\f pm]$.}
\endcases\endaligned$$\endpro
Proof. Suppose $p=2mk+r$ with $k\in\Bbb Z$ and
$r\in\{0,1,\ldots,2m-1\}$. Then $[\f pm]=[\f rm]$ and so $2\mid [\f
pm]$ if and only if $r<m$. Hence
$$\aligned &\Big[\f
p{2m}\Big]+\Big[\f{(m-1)p}{2m}\Big]\\&=k+k(m-1)+\Big[\f{(m-1)r}{2m}\Big]=\f{p-r}2+
\Big[\f{(m-1)r}{2m}\Big]
\\&=\f{p-1}2+\Big[\f{(m-1)r+m(1-r)}{2m}\Big]=\f{p-1}2+\Big[\f{m-r}{2m}\Big]
\\&=\cases \f{p-1}2&\t{if $2\mid [\f pm]$,}
\\\f{p-3}2&\t{if $2\nmid [\f pm]$.}\endcases
\endaligned$$
Thus, if $2\mid [\f pm]$, using Lemma 2.1 and the above we
get
$$\align P_{[\f pm]}(\sqrt t)&=P_{2[\f p{2m}]}(\sqrt
t)=\sum_{k=0}^{[\f p{2m}]}\b{[\f p{2m}]}k\b{-\f 12-[\f
p{2m}]}k(1-t)^k
\\&\e \sum_{k=0}^{[\f p{2m}]}\b{[\f
p{2m}]}k\b{\f {p-1}2-[\f p{2m}]}k(1-t)^k
\\&=\sum_{k=0}^{[\f p{2m}]}\b{[\f
p{2m}]}k\b{[\f{(m-1)p}{2m}]}k(1-t)^k\mod p.\endalign$$ If $2\nmid
[\f pm]$, using Lemma 2.1 and the above we get
$$\align P_{[\f pm]}(\sqrt t)/\sqrt t&=P_{2[\f p{2m}]+1}(\sqrt
t)/\sqrt t=\sum_{k=0}^{[\f p{2m}]}\b{[\f p{2m}]}k\b{-\f 32-[\f
p{2m}]}k(1-t)^k
\\&\e \sum_{k=0}^{[\f p{2m}]}\b{[\f
p{2m}]}k\b{\f {p-3}2-[\f p{2m}]}k(1-t)^k
\\&=\sum_{k=0}^{[\f p{2m}]}\b{[\f
p{2m}]}k\b{[\f{(m-1)p}{2m}]}k(1-t)^k\mod p.\endalign$$ By Lemma 2.1,
we also have
$$\aligned
&t^{[\f p{2m}]}\sum_{k=0}^{[\f p{2m}]}\b{[\f p{m}]}{2k}\b{\f{p-1}2}k
\Big(\f 1t-1\Big)^k
\\&\e (\sqrt t)^{2[\f p{2m}]}\sum_{k=0}^{[\f p{2m}]}\b{[\f
p{m}]}{2k}\b{-\f{1}2}k \Big(\f 1t-1\Big)^k
\\&=\cases (\sqrt t)^{[\f pm]}\sum_{k=0}^{[\f p{2m}]}\b{[\f
p{m}]}{2k}\b{-\f{1}2}k (\f 1t-1)^k=P_{[\f pm]}(\sqrt t)\mod p &\t{if
$2\mid [\f pm]$,}
\\(\sqrt t)^{[\f pm]-1}\sum_{k=0}^{[\f p{2m}]}\b{[\f
p{m}]}{2k}\b{-\f{1}2}k (\f 1t-1)^k=\f{P_{[\f pm]}(\sqrt t)}{\sqrt
t}\mod p &\t{if $2\nmid [\f pm]$.}\endcases\endaligned$$ This
completes the proof.

\pro{Lemma 2.2} Let $p$ be an odd prime and
$m\in\{1,2,\ldots,\f{p-1}2\}$. Then $P_{p-1-m}(x)\e P_m(x)$.
\endpro
Proof. Since $m<\f p2$ we have $p-1-m\ge m$. Note that
$\b{-t}k=(-1)^k\b{t+k-1}k$. By using (1.3) we have
$$\align P_{p-1-m}(x)&=\sum_{k=0}^{p-1-m}
\b{p-1-m}k\b{p-1-m+k}k\Ls{x-1}2^k \\&\e \sum_{k=0}^{p-1-m}
\b{-1-m}k\b{-1-m+k}k\Ls{x-1}2^k
\\&=\sum_{k=0}^{p-1-m}(-1)^k\b{m+k}k\cdot (-1)^k\b mk\Ls{x-1}2^k
\\&=\sum_{k=0}^m\b mk\b{m+k}k\Ls{x-1}2^k
=P_m(x)\mod p.\endalign$$ This proves the lemma.

 \pro{Lemma 2.3} Let
$p$ be an odd prime and $m\in\{0,1,\ldots,p-1\}$. Then
  $$P_{p+m}(x)\e x^pP_m(x)\mod p.$$
 \endpro
 Proof. If $a_i,b_i\in\{0,1,\ldots,p-1\}$, then Lucas theorem
 asserts that
 $$\b{a_0+a_1p+\cdots+a_np^n}{b_0+b_1p+\cdots+b_np^n}
 \e \b{a_0}{b_0}\b{a_1}{b_1}\cdots\b{a_n}{b_n}\mod p.$$ Thus,
for $0\le r\le m<p$ we have
$$\b{2p+m+r}{p+r}\e\cases \b 21\b{m+r}r\mod p&\t{if $r<p-m$,}
\\\b {3p}{p+r}\e \b 31\b 0r=0\e \b 21\b{m+r}r\mod p&\t{if $r=p-m$.}
\\\b 31\b{m+r-p}r=0\e \b 21\b{m+r}r\mod p&\t{if $r> p-m$.}
\endcases$$
Hence, using (1.3) we see that
 $$\align P_{p+m}(x)&=\sum_{k=0}^{p+m}\b{p+m}k
 \b{p+m+k}k\Ls{x-1}2^k
\\&=\sum_{k=0}^{p-1}\b{p+m}k
 \b{p+m+k}k\Ls{x-1}2^k\\&\q+\sum_{k=p}^{p+m}\b{p+m}k
 \b{p+m+k}k\Ls{x-1}2^k
 \\&\e \sum_{k=0}^{p-1}\b{m}k
 \b{m+k}k\Ls{x-1}2^k\\&\q+\sum_{r=0}^m\b{p+m}{p+r}
 \b{2p+m+r}{p+r}\Ls{x-1}2^{p+r}
 \\&\e \sum_{k=0}^m\b{m}k
 \b{m+k}k\Ls{x-1}2^k+\sum_{r=0}^m\b mr
 \cdot 2\b{m+r}r\Ls{x-1}2^{p+r}
\\&\e (1+(x-1)^p)\sum_{k=0}^m\b{m}k
 \b{m+k}k\Ls{x-1}2^k\\&\e x^pP_m(x)\mod p.\endalign$$
 Thus the lemma is proved.
\pro{Theorem 2.2} Let $p$ be an odd prime and $a,t\in\Bbb Z_p$ with
$t\not\e 0\mod p$. Then
$$\aligned P_{2\ap}(\sqrt{t})&\e\sum_{k=0}^{p-1}\b ak\b{-\f 12-a}k(1-t)^k
\\&\e t^{\ap}\sum_{k=0}^{p-1}\b ak\b{a-\f 12}k\Big(1-\f 1t\Big)^k
\\&\e
\Ls tp t^{[\f{2\ap}p]}\sum_{k=0}^{p-1}\b{-1-a}k\b{a-\f 12}k(1-t)^k
\\&\e \Ls tpt^{[\f {2\ap}p]-\ap}\sum_{k=0}^{p-1}\b{-1-a}k\b{-\f
12-a}k\Big(1-\f 1t\Big)^k \mod p. \endaligned$$
\endpro
Proof.  For $\beta\in\Bbb Z_p$ we have
$$\align& \sum_{k=0}^{p-1}\b ak\b{a+\beta}k\Ls{t-1}t^k
\\&\e \sum_{k=0}^{p-1}\b{\ap}k\b{\ap+\beta}k\Ls{t-1}t^k
\\&=\sum_{k=0}^{\ap}\b{\ap}k\b{\ap+\beta}k\Ls{t-1}{t}^k
\\&=\sum_{k=0}^{\ap}\b{\ap}{\ap-k}\b{\ap+\beta}{\ap-k}\Ls{t-1}t^{\ap-k}
\\&=\f 1{{(2t)}^{\ap}}\sum_{k=0}^{\ap}\b{\ap}k\b{\ap+\beta}{\ap-k}
((2t-1)+1)^k((2t-1)-1)^{\ap-k}
\\&=\f 1{t^{\ap}}P_{\ap}^{(0,\beta)}(2t-1)\mod p.\endalign$$
Thus, by (2.2) and the fact $\ag{-1-a}=p-1-\ap$ we get
 $$\align\sum_{k=0}^{p-1}\b ak\b{a-\f 12}k\Ls
{t-1}t^k \e \f 1{t^{\ap}}P_{\ap}^{(0,-\f 12)}(2t-1) =\f
1{t^{\ap}}P_{2\ap}(\sqrt t)\mod p\endalign$$ and
$$\align&\sum_{k=0}^{p-1}\b {-1-a}k\b{-\f 12-a}k\Ls
{t-1}t^k \\&\e t^{-\ag{-1-a}}P_{\ag{-1-a}}^{(0,\f 12)}(2t-1) =
t^{-\ag{-1-a}}P_{2\ag{-1-a}+1}(\sqrt t)/\sqrt t
\\&\e t^{\ap}P_{2(p-1-\ap)+1}(\sqrt t)/\sqrt t\mod p.\endalign$$

Using (2.3) and (2.2) we see that
$$\align &\sum_{k=0}^{p-1}\b ak\b{-\f 12-a}k(1-t)^k
\\&\e \sum_{k=0}^{p-1}\b {\ap}k\b{-\f 12-\ap}k(1-t)^k
= \sum_{k=0}^{\ap}\b {\ap}k\b{-\f 12-\ap}k(1-t)^k
\\&=P_{\ap}^{(0,-\f 12)}(2t-1)=P_{2\ap}(\sqrt{t})\mod p
\endalign$$
and
$$\align &\sum_{k=0}^{p-1}\b{-1-a}k\b{a-\f 12}k(1-t)^k
\\&\e \sum_{k=0}^{p-1}\b{p-1-\ap}k\b{1+a-\f 12-1}k(1-t)^k
\\&\e \sum_{k=0}^{p-1-\ap}\b{p-1-\ap}k\b{-(p-1-\ap)-\f 12-1}k(1-t)^k
\\&=P_{p-1-\ap}^{(0,\f 12)}(2t-1)=P_{2(p-1-\ap)+1}(\sqrt t)/\sqrt
t\mod p.\endalign$$ To complete the proof, using Lemmas 2.2 and 2.3
we note that
$$\aligned &P_{2(p-1-\ap)+1}(\sqrt t)/\sqrt
t\\&\e\cases \f{P_{p-1-2\ap}(\sqrt t)(\sqrt t)^p}{\sqrt t} \e \sls
tpP_{2\ap}(\sqrt t)\mod p\q\t{if $\ap<\f
p2$,}\\\f{P_{p-1-(2\ap-p)}(\sqrt t)}{\sqrt t} \e \f{P_{2\ap-p}(\sqrt
t)}{\sqrt t}\e \f{P_{2\ap}(\sqrt t)}{(\sqrt t)^{p+1}} \\\q\e \sls
tp\f 1t P_{2\ap}(\sqrt t)\mod p\qq\qq\qq\q\ \;\t{if $\ap>\f p2$.}
\endcases\endaligned$$

\pro{Theorem 2.3} Let $p$ be an odd prime and $a,t\in\Bbb Z_p$ with
$t\not\e 0,1\mod p$. Then
$$\align\sum_{k=0}^{p-1}\b{a}k^2t^k&\e t^{\ap}
\sum_{k=0}^{p-1}\f{\b{a}k^2}{t^k}\e (t-1)^{\ap}P_{\ap}
\Ls{t+1}{t-1}
\\&\e (t-1)^{\ap}\sum_{k=0}^{p-1}\b ak\b{-1-a}k\f 1{(1-t)^k}
\\&\e (t-1)^{2\ap}\sum_{k=0}^{p-1}\b{-1-a}k^2t^k
\mod p.\endalign$$
\endpro
Proof. It is clear that
$$\align\sum_{k=0}^{p-1}\b{a}k^2t^k
&\e \sum_{k=0}^{\ap}\b{\ap}k^2t^k
=\sum_{s=0}^{\ap}\b{\ap}{\ap-s}^2t^{\ap-s}
\\&=t^{\ap}\sum_{k=0}^{\ap}\b{\ap}k^2t^{-k}
\e t^{\ap} \sum_{k=0}^{p-1}\f{\b{a}k^2}{t^k}\mod p.\endalign$$ Using
(2.1) we see that for $x\not=1$,
$$\align P_{\ap}(x)&=P_{\ap}^{(0,0)}(x)=\Ls{x-1}2^{\ap}\sum_{k=0}^{\ap}
\b {\ap}k^2\Ls{x+1}{x-1}^k\\& =\Ls{x-1}2^{\ap}\sum_{k=0}^{p-1} \b
{\ap}k^2\Ls{x+1}{x-1}^k
\\&\e \Ls{x-1}2^{\ap}\sum_{k=0}^{p-1} \b
{a}k^2\Ls{x+1}{x-1}^k\mod p.\endalign$$ Set $x=\f{t+1}{t-1}$. Then
$t=\f{x+1}{x-1}$ and $\f{x-1}2=\f 1{t-1}$. Now substituting $x$ with
$\f{t+1}{t-1}$ in the above congruence we obtain
$$\sum_{k=0}^{p-1}\b{a}k^2t^k\e (t-1)^{\ap}P_{\ap}
\Ls{t+1}{t-1}\mod p.$$
 Clearly $\ag{-1-a}=p-1-\ap$. Thus, using Lemma 2.2 and the above
we see that
$$\align &\sum_{k=0}^{p-1}\b {-1-a}k^2t^k
\\&\e (t-1)^{\ag{-1-a}}P_{\ag{-1-a}}
\Ls{t+1}{t-1}=(t-1)^{p-1-\ap}P_{p-1-\ap}\Ls{t+1}{t-1}
\\&\e (t-1)^{-\ap}P_{\ap}\Ls{t+1}{t-1}
\e (t-1)^{-2\ap}\sum_{k=0}^{p-1}\b ak^2t^k \mod p.\endalign$$ To
complete the proof, using (2.3) we note that
$$\align P_{\ap}\Ls{t+1}{t-1}&=\sum_{k=0}^{\ap}\b{\ap}k\b{-1-\ap}k\f
1{(1-t)^k}\\&=\sum_{k=0}^{p-1}\b{\ap}k\b{-1-\ap}k\f 1{(1-t)^k}
\\&\e \sum_{k=0}^{p-1}\b ak\b {-1-a}k\f
1{(1-t)^k}\mod p.\endalign$$

\subheading{3. Congruences for $\sum_{x=0}^{p-1}\sls{x^3+mx+n}p\mod
p$} \pro{Theorem 3.1} Let $p$ be an odd prime and let $t\in\Bbb
Z_p$. Then
$$\aligned&\sum_{k=0}^{p-1}\b{-\f 1{12}}k\b{-\f 5{12}}k(1-t)^k
\\&\e t^{\ag{-\f 1{12}}}\sum_{k=0}^{p-1}\b{-\f 1{12}}k\b{-\f 7{12}}k\Big(1-\f 1t\Big)^k
\\&\e \cases P_{[\f p6]}(\sqrt t)\mod p&\t{if $p\e 1\mod 4$,}
\\\sls tpP_{[\f p6]}(\sqrt{t})\sqrt{t}\mod p&\t{if $p\e 3\mod
4$.}\endcases\endaligned$$
\endpro
Proof. Taking $a=-\f 1{12}$ in Theorem 2.2 and then applying Lemmas
2.2-2.3 we see that
$$\aligned&\sum_{k=0}^{p-1}\b{-\f 1{12}}k\b{-\f 5{12}}k(1-t)^k
\\&\e t^{\ag{-\f 1{12}}} \sum_{k=0}^{p-1}\b{-\f 1{12}}k\b{-\f
7{12}}k\Ls {t-1}t^k \e P_{2\ag{-\f 1{12}}}(\sqrt t)
\\&=\cases P_{\f{p-1}6}(\sqrt{t})\mod p&\t{if $12\mid p-1$,}
\\P_{\f{5p-1}6}(\sqrt{t})\e P_{\f{p-5}6}(\sqrt{t})\mod p&\t{if
$12\mid p-5$,}
\\P_{\f{7p-1}6}(\sqrt{t})\e
(\sqrt{t})^pP_{\f{p-1}6}(\sqrt{t})\mod p& \t{if $12\mid p-7$,}
\\P_{\f{11p-1}6}(\sqrt{t})\e
(\sqrt{t})^pP_{\f{5p-1}6} (\sqrt{t})\e(\sqrt{t})^p
P_{\f{p-5}6}(\sqrt{t})\mod p&\t{if $12\mid p-11$.}
\endcases\endaligned$$
To see the result, we note that $(\sqrt{t})^p=t^{\f{p-1}2}\sqrt{t}\e
\sls {t}p\sqrt{t}\mod p$.

\pro{Lemma 3.1 ([S5, Theorem 2.2])} Let $p>3$ be a prime and
$m,n\in\Bbb Z_p$ with $m\not\e 0\mod p$. Then
$$\sum_{x=0}^{p-1}\Ls{x^3+mx+n}p\e \cases
-(-3m)^{\f{p-1}4}P_{[\f p6]}\ls{3n\sqrt{-3m}}{2m^2}\mod p&\t{if $p\e
1\mod 4$,}
\\-\f{(-3m)^{\f{p+1}4}}{\sqrt{-3m}}P_{[\f p6]}\ls{3n\sqrt{-3m}}{2m^2}\mod
p&\t{if $p\e 3\mod 4$.}
\endcases$$
\endpro
\pro{Theorem 3.2} Let $p>3$ be a prime and $m,n\in\Bbb Z_p$ with
$mn\not\e 0\mod p$. Then
$$\aligned &\sum_{x=0}^{p-1}\Ls{x^3+mx+n}p
\\&\e\cases  -(-3m)^{\f{p-1}4}
\sum_{k=0}^{p-1}\b{-\f 1{12}}k\b{-\f 5{12}}k
(\f{4m^3+27n^2}{4m^3})^k\mod p&\t{if $4\mid p-1$,}
\\\f{2m}{9n}\sls{-3m}p(-3m)^{\f{p+1}4}
\sum_{k=0}^{p-1}\b{-\f 1{12}}k\b{-\f 5{12}}k
(\f{4m^3+27n^2}{4m^3})^k\mod p&\t{if $4\mid p-3$.}
\endcases
\\&\e\cases  (-1)^{\f{p+1}2}\sls n2^{\f{p-1}6}
\sum_{k=0}^{p-1}\b{-\f 1{12}}k\b{-\f 7{12}}k
\sls{4m^3+27n^2}{27n^2}^k\mod p&\t{if $3\mid p-1$,}
\\ (-1)^{\f{p+1}2}\f 3m\sls 2n^{\f{p-5}6}
\sum_{k=0}^{p-1}\b{-\f 1{12}}k\b{-\f 7{12}}k
\sls{4m^3+27n^2}{27n^2}^k\mod p&\t{if $3\mid p-2$.}\endcases
\endaligned$$\endpro
Proof. Set $t=-\f{27n^2}{4m^3}$. Then $1-t=\f{4m^3+27n^2}{4m^3}$.
Since
$$t^{\ag{-\f 1{12}}}=\cases
(-\f{27n^2}{4m^3})^{\f{p-1}{12}}=(-\f 3m)^{\f{p-1}4}\sls
n2^{\f{p-1}6}&\t{if $p\e 1\mod{12}$,}
\\(-\f{27n^2}{4m^3})^{\f{5p-1}{12}}\e\sls{-3m}p(-\f
m3)^{\f{p-5}4}\sls 2n^{\f{p-5}6}\mod p&\t{if $p\e 5\mod{12}$,}
\\(-\f{27n^2}{4m^3})^{\f{7p-1}{12}}\e\sls{-3m}p(-\f
3m)^{\f{p+5}4}\sls n2^{\f{p+5}6}\mod p&\t{if $p\e 7\mod{12}$,}
\\(-\f{27n^2}{4m^3})^{\f{11p-1}{12}}\e(-\f
m3)^{\f{p-11}4}\sls 2n^{\f{p-11}6}\mod p&\t{if $p\e 11\mod{12}$,}
\endcases$$ using Theorem 3.1 and Lemma 3.1 we deduce the result.

\pro{Corollary 3.1} Let $p$ be a prime with $p\not=2,3,19$. Then
$$\aligned&\sum_{k=0}^{p-1}\b{-\f 1{12}}k\b{-\f 7{12}}k\f 1{513^k}
\\&\e
\cases 0\mod p&\t{if $\sls p{19}=-1$,}
\\-19^{-\f{p-1}3}\sls{-2}p\sls u{19}u\mod p&\t{if
$\sls p{19}=\sls p3=1$ and so $4p=u^2+19v^2$.}
\\\f 83\cdot 19^{\f {p-2}3}\sls{-2}p\sls
u{19}u\mod p&\t{if $\sls p{19}=-\sls p3=1$ and so $4p=u^2+19v^2$.}
\endcases\endaligned$$
\endpro
Proof. From [RPR], [JM] and [PV] we know that
$$\aligned \sum_{x=0}^{p-1}\Ls{x^3-8\cdot 19x+2\cdot 19^2}p
=\cases \sls 2p\sls u{19}u&\t{if $\sls p{19}=1$ and
$4p=u^2+19v^2$,}\\0&\t{if $\sls p{19}=-1$.}
\endcases\endaligned\tag 3.1$$
Thus, taking $m=-8\cdot 19$ and $n=2\cdot 19^2$ in Theorem 3.2 we
deduce the result.
\par For a simpler congruence for $u\mod p$ with $4p=u^2+19v^2$,
see [LH, p.269].

\pro{Theorem 3.3} Let $p>3$ be a prime and $m,n\in\Bbb Z_p$ with
$mn\not\e 0\mod p$. Then
$$\aligned &\sum_{x=0}^{p-1}\Ls{x^3+mx+n}p
\\&\e\cases  -(-3m)^{\f{p-1}4}
\sum_{k=0}^{[p/12]}\b{[\f p{12}]}k\b{[\f {5p}{12}]}k
(\f{4m^3+27n^2}{4m^3})^k\mod p&\t{if $4\mid p-1$,}
\\-\f{3n}{2m^2}(-3m)^{\f{p+1}4}
\sum_{k=0}^{[p/12]}\b{[\f p{12}]}k\b{[\f {5p}{12}]}k
(\f{4m^3+27n^2}{4m^3})^k\mod p&\t{if $4\mid p-3$}
\endcases
\\&\e\cases (-1)^{\f{p+1}2}\sls n2^{\f{p-1}6}\sum_{k=0}^{[p/12]}\b{[\f
p{6}]}{2k}\b{\f{p-1}2}k(-\f{4m^3+27n^2}{27n^2})^k\mod p&\t{if $3\mid
p-1$,}
\\(-1)^{\f{p+1}2}\f m3\sls n2^{\f{p-5}6}\sum_{k=0}^{[p/12]}\b{[\f
p{6}]}{2k}\b{\f{p-1}2}k(-\f{4m^3+27n^2}{27n^2})^k\mod p&\t{if $3\mid
p-2$.}\endcases
\endaligned$$
\endpro
Proof. Putting $m=6$ and $t=-\f{27n^2}{4m^3}$ in Theorem 2.1 we see
that
$$\aligned&\sum_{k=0}^{[p/12]}\b{[\f p{12}]}k\b{[\f {5p}{12}]}k
\Big(\f{4m^3+27n^2}{4m^3}\Big)^k
\\&\e \Big(-\f{27n^2}{4m^3}\Big)^{[\f p{12}]}\sum_{k=0}^{[p/12]}\b{[\f
p{6}]}{2k}\b{\f{p-1}2}k\Big(-\f{4m^3+27n^2}{27n^2}\Big)^k
\\&\e\cases P_{[\f
p6]}(\f{3n\sqrt{-3m}}{2m^2})\mod p&\t{if $p\e 1,5\mod {12}$,}
\\P_{[\f
p6]}(\f{3n\sqrt{-3m}}{2m^2})/(\f{3n\sqrt{-3m}}{2m^2})\mod p&\t{if
$p\e 7,11\mod {12}$.}\endcases\endaligned$$ Now applying Lemma 3.1
we deduce the result.

\par\q
\par Let $p>3$ be a prime and $S_p(m,n)=\sum_{x=0}^{p-1}\sls{x^3+mx+n}p$. It is known that
(see for example [S5, S6], [R1,R2], [RPR], [JM], [PV] and [W])
$$\aligned S_p(-11,14)
=\cases (-1)^{\f{p+3}4}2a&\t{if $4\mid p-1$, $p=a^2+b^2$ and $4\mid
a-1$,}
\\0&\t{if $p\e 3\mod 4$,}
\endcases\endaligned\tag 3.2$$
$$\aligned S_p(-30,56)
= \cases (-1)^{\f{p+7}8}\sls 3p2c&\t{if $p\e 1\mod 8$, $p=c^2+2d^2$
and $4\mid c-1$,}
\\(-1)^{\f{p-3}8}\sls 3p2c&\t{if $p\e 3\mod 8$,
$p=c^2+2d^2$ and $4\mid c-1$,}
\\0&\t{if $p\e 5,7\mod 8$,}
\endcases\endaligned\tag 3.3$$
$$\aligned S_p(-15,22)
=\cases -2A&\t{if $3\mid p-1$, $p=A^2+3B^2$ and $3\mid A-1$,}
\\0&\t{if $p\e 2\mod 3$,}
\endcases\endaligned\tag 3.4$$
$$S_p(-120,506)=\cases \sls 2pL&\t{if $3\mid
p-1$, $4p=L^2+27M^2$ and $3\mid L-1$,}\\0&\t{if $p\e 2\mod 3$,}
\endcases\tag 3.5$$
$$\aligned S_p(-35,98)
=\cases(-1)^{\f{p+1}2}2\sls C7C&\t{if $p\e 1,2,4\mod 7$ and
 $p=C^2+7D^2$,}\\0&\t{if $p\e 3,5,6\mod
7$,}\endcases\endaligned\tag 3.6$$
$$\aligned S_p(-595,5586)=\cases
 (-1)^{\f{p+1}2}2C\sls C7&\t{if $p=C^2+7D^2\e 1,2,4\mod 7$,}
\\0&\t{if $p\e 3,5,6\mod 7$,}
\endcases\endaligned\tag 3.7$$
$$S_p(-96\cdot 11,112\cdot 11^2)=\cases \sls 2p\sls
u{11}u&\t{if $\sls p{11}=1$ and $4p=u^2+11v^2$,}\\0&\t{if $\sls
p{11}=-1$.}\endcases\tag 3.8$$
$$ S_p(-80\cdot 43,42\cdot 43^2)
=\cases \sls 2p\sls u{43}u&\t{if $\sls p{43}=1$ and
$4p=u^2+43v^2$,}\\0&\t{if $\sls p{43}=-1$,}
\endcases\tag 3.9$$
$$S_p(-440\cdot 67,434\cdot 67^2) =\cases
\sls 2p\sls u{67}u&\t{if $\sls p{67}=1$ and
$4p=u^2+67v^2$,}\\0&\t{if $\sls p{67}=-1$,}
\endcases\tag 3.10$$
$$\aligned&S_p(-80\cdot 23\cdot 29\cdot
163,14\cdot 11\cdot 19\cdot 127\cdot 163^2)\\&=\cases \sls 2p\sls
u{163}u&\t{if $\sls p{163}=1$ and $4p=u^2+163v^2$,}\\0&\t{if $\sls
p{163}=-1$.}
\endcases
\endaligned \tag 3.11$$
Using (3.2)-(3.11) and Theorems 3.2 and 3.3 we may deduce the
results similar to Corollary 3.1.

\subheading {4. Congruences for $\sum_{k=0}^{[p/6]}\b{[\f
p3]}k\b{[\f p6]}k(1-t)^k\mod p$}

 \pro{Theorem 4.1} Let $p$ be a prime with
$p\e 1,4\mod 5$. Then
$$\aligned &\sum_{k=0}^{[p/6]}\b{[\f p3]}k\b{[\f p6]}k(-4)^k
\\&\e\cases 2x\mod p&\t{if $p\e 1,4\mod{15}$ and so
$p=x^2+15y^2$ with $3\mid x-1$,}
\\0\mod p&\t{if $p\e 11,14\mod {15}$.}\endcases\endaligned$$
\endpro
Proof. Taking $m=3$ and $t=5$ in Theorem 2.1 we get
$$\sum_{k=0}^{[p/6]}\b{[\f p3]}k\b{[\f p6]}k(-4)^k\e
\cases P_{[\f p3]}(\sqrt 5)\mod p&\t{if $p\e 1\mod 3$,}
\\P_{\f{p-2}3}(\sqrt 5)/\sqrt 5\mod p&\t{if $p\e 2\mod 3$.}\endcases$$
From [S6, Theorem 4.6] we know that
$$P_{[\f p3]}(\sqrt 5)\e \cases 2x\mod p&\t{if $p=x^2+15y^2\e 1,4\mod {15}$ and $3\mid
x-1$,}\\0\mod p&\t{if $p\e 11,14\mod{15}$.}
\endcases$$
Thus the result follows.
 \pro{Conjecture 4.1} Let $p>5$ be a prime.
Then
$$\aligned\sum_{k=0}^{p-1}\b{-\f 13}k\b{-\f 16}k(-4)^k
&\e \Ls 5p5^{\f{1-\sls p3}2}\sum_{k=0}^{p-1}\b{-\f 23}k\b{-\f
56}k(-4)^k\\& \e\cases \sls x3(2x-\f p{2x})\mod {p^2}&\t{if
$p=x^2+15y^2$,}
\\ -\sls x3(10x-\f p{2x})\mod {p^2}&\t{if $p=5x^2+3y^2$,}
\\0\mod p&\t{if $p\e 17,23\mod{30}$}
\endcases\endaligned$$ and so
$$2x\Ls x3\e
\sum_{k=0}^{(p-5)/6}\b{\f{p-2}3}k\b{\f{p-5}6}k(-4)^k\mod
p\qtq{for}p=5x^2+3y^2.$$
\endpro
\pro{Theorem 4.2} Let $p$ be a prime such that $p\e 1,7\mod 8$. Then
$$\aligned
&2^{[\f p3]}\sum_{k=0}^{[p/6]}\b{[\f p3]}{2k}\b{\f {p-1}2}k\e
\sum_{k=0}^{[p/6]}\b{[\f p3]}k\b{[\f p6]}k\f 1{2^k}
\\&\e\cases 2x\mod p&\t{if $p\e 1,7\mod{24}$ and so
$p=x^2+6y^2$ with $3\mid x-1$,}\\0\mod p&\t{if $p\e 17,23\mod{24}$.}
\endcases\endaligned$$\endpro
Proof. Taking $m=3$ and $t=\f 12$ in Theorem 2.1 we see that
$$\aligned&\sum_{k=0}^{[p/6]}\b{[\f p3]}k\b{[\f p6]}k\f 1{2^k}
\e \f 1{2^{[p/6]}}\sum_{k=0}^{[p/6]}\b{[\f p3]}{2k} \b{\f{p-1}2}k
\\&
\e\cases P_{\f{p-1}3}(\f 1{\sqrt 2})\mod p&\t{if $p\e 1\mod 3$,}
\\\sqrt 2P_{\f{p-2}3}(\f 1{\sqrt 2})\mod p&\t{if $p\e 2\mod
3$.}
\endcases\endaligned$$
From [S6, Theorem 4.5] we have
$$P_{[\f p3]}\Ls 1{\sqrt 2}\e \cases
2x\sls x3\mod p&\t{if $p=x^2+6y^2\e 1,7\mod {24}$,}
 \\0\mod p&\t{if
$p\e 17,23\mod{24}$}\endcases$$ Observe that
$2^{-\f{p-1}6}=2^{\f{p-1}3-\f{p-1}2}\e 2^{\f{p-1}3}\mod p$ for $p\e
1,7\mod {24}$. From the above we deduce the result.

\pro{Conjecture 4.2} Let $p>3$ be a prime. Then
$$\aligned&\sum_{k=0}^{p-1}\f{\b{-\f 13}k\b{-\f 16}k}{2^k}
\e \Ls 2p2^{\f{\sls p3-1}2}\sum_{k=0}^{p-1}\f{\b{-\f 23}k\b{-\f
56}k}{2^k}
\\&\e\cases \sls x3(2x-\f p{2x})\mod {p^2}&\t{if $p=x^2+6y^2\e 1,7\mod
{24}$,}\\\sls x3(2x-\f p{4x})\mod {p^2}&\t{if $p=2x^2+3y^2\e
5,11\mod {24}$,}
\\0\mod p&\t{if $p\e 13,19\mod {24}$}\endcases\endaligned$$
and so
$$x\Ls x3\e -\f 14
\sum_{k=0}^{(p-5)/6}\b{\f{p-2}3}k\b{\f{p-5}6}k\f 1{2^k}\mod
p\qtq{for}p=2x^2+3y^2.$$
\endpro
\pro{Theorem 4.3} Let $p$ be an odd prime such that $\sls{17}p=1$.
Then
$$\aligned&\sum_{k=0}^{[p/6]}\b{[\f p3]}k\b{[\f p6]}k\f 1{(-16)^k}
\\&\e\cases x\mod p&\t{if $p\e 1\mod 3$ and so $4p=x^2+51y^2$ with $3\mid x-2$,}
\\0\mod p&\t{if $p\e 2\mod 3$.}
\endcases\endaligned$$\endpro
Proof. Taking $m=3$ and $t=\f {17}{16}$ in Theorem 2.1 we see that
$$\sum_{k=0}^{[p/6]}\b{[\f p3]}k\b{[\f p6]}k\f 1{(-16)^k}
\e\cases P_{[\f p3]}(\f{\sqrt{17}}4)\mod p&\t{if $p\e 1\mod 3$,}
\\\f 4{\sqrt{17}}P_{[\f p3]}(\f{\sqrt{17}}4)
\mod p&\t{if $p\e 2\mod 3$.}\endcases$$ On the other hand, by [S6,
Theorem 4.8],
$$P_{[\f p3]}\Ls{\sqrt{17}}4\e \cases  -\sls x3x\mod p&\t{if $p\e 1\mod 3$ and so
$4p=x^2+51y^2$,}
\\0\mod p&\t{if $p\e 2\mod 3$.}\endcases$$ Thus the theorem is
proved.

\pro{Conjecture 4.3} Let $p>3$ be a prime. Then
$$\aligned\sum_{k=0}^{p-1}\f{\b{-\f 13}k\b{-\f 16}k}{(-16)^k}
&\e\Ls{17}p\Ls{17}{16}^{\f{1-\sls p3}2}\sum_{k=0}^{p-1}\f{\b{-\f
23}k\b{-\f 56}k}{(-16)^k} \\&\e\cases -\sls x3(x-\f px)\mod
{p^2}&\t{if $4p=x^2+51y^2$,}
\\ \f 14\sls x3(17x-\f px)\mod{p^2}&\t{if $4p=17x^2+3y^2$,}
\\0\mod p&\t{if $\sls p3=-\sls p{17}=1$}\endcases\endaligned$$
and so
$$x\Ls x3\e -\f 14
\sum_{k=0}^{(p-5)/6}\b{\f{p-2}3}k\b{\f{p-5}6}k\f 1{(-16)^k}\mod
p\qtq{for}4p=17x^2+3y^2.$$

\endpro

\par Using the theorems in Section 4 in [S6] and Theorem 2.1 one can
similarly deduce the following results.
 \pro{Theorem 4.4} Let $p$ be
an odd prime such that $\sls{41}p=1$. Then
$$\aligned&\sum_{k=0}^{[p/6]}\b{[\f p3]}k\b{[\f p6]}k\f 1{(-1024)^k}
\\&\e\cases x\mod p&\t{if $p\e 1\mod 3$ and so $4p=x^2+123y^2$ with $3\mid x-2$,}
\\0\mod p&\t{if $p\e 2\mod 3$.}
\endcases\endaligned$$\endpro

\pro{Conjecture 4.4} Let $p>3$ be a prime. Then
$$\aligned\sum_{k=0}^{p-1}\f{\b{-\f 13}k\b{-\f 16}k}{(-1024)^k}
&\e \Ls{41}p\Ls{1025}{1024}^{\f{1-\sls p3}2}
\sum_{k=0}^{p-1}\f{\b{-\f 23}k\b{-\f 56}k}{(-1024)^k}
 \\&\e\cases -\sls x3(x-\f px)\mod {p^2}&\t{if $4p=x^2+123y^2$,}
\\ \f 5{32}\sls x3(41x-\f px)\mod{p^2}&\t{if $4p=41x^2+3y^2$,}
\\0\mod p&\t{if $\sls p3=-\sls p{41}=1$}\endcases\endaligned$$

and so
$$x\Ls x3\e
-\f 5{32}\sum_{k=0}^{(p-5)/6}\b{\f{p-2}3}k\b{\f{p-5}6}k\f
1{(-1024)^k}\mod p\qtq{for}4p=41x^2+3y^2.$$

\endpro

\pro{Theorem 4.5} Let $p$ be an odd prime such that $\sls{89}p=1$.
Then
$$\aligned&\sum_{k=0}^{[p/6]}\b{[\f p3]}k\b{[\f p6]}k\f 1{(-250000)^k}
\\&\e\cases x\mod p&\t{if $p\e 1\mod 3$ and so $4p=x^2+267y^2$ with $3\mid x-2$,}
\\0\mod p&\t{if $p\e 2\mod 3$.}
\endcases\endaligned$$\endpro

\pro{Conjecture 4.5} Let $p>3$ be a prime. Then
$$\aligned\sum_{k=0}^{p-1}\f{\b{-\f 13}k\b{-\f 16}k}{(-250000)^k}
&\e \Ls {89}p\Ls{250001}{250000}^{\f{1-\sls
p3}2}\sum_{k=0}^{p-1}\f{\b{-\f 23}k\b{-\f 56}k}{(-250000)^k}
 \\&\e\cases -\sls x3(x-\f px)\mod {p^2}&\t{if $4p=x^2+267y^2$,}
\\ \f{53}{500}\sls x3(89x-\f px)\mod{p^2}&\t{if $4p=89x^2+3y^2$,}
\\0\mod p&\t{if $\sls p3=-\sls p{89}=1$}\endcases\endaligned$$\endpro
and so
$$x\Ls x3\e -\f{53}{500}
\sum_{k=0}^{(p-5)/6}\b{\f{p-2}3}k\b{\f{p-5}6}k\f 1{(-250000)^k}\mod
p\ \t{for}\ 4p=89x^2+3y^2.$$

\pro{Theorem 4.6}  Let $p$ be a prime such that $p\e 1,4\mod 5$.
Then
$$\aligned&\sum_{k=0}^{[p/6]}\b{[\f p3]}k\b{[\f p6]}k\f 1{(-80)^k}
\\&\e \cases x\mod
p&\t{if $p\e 1,4\mod {15}$ and so $4p=x^2+75y^2$ with $3\mid x-2$,}
 \\0\mod p&\t{if
$p\e 11,14\mod{15}$}\endcases\endaligned$$
\endpro

\pro{Conjecture 4.6}  Let $p>5$ be a prime. Then
$$\aligned&\sum_{k=0}^{p-1}\b{-\f 13}k\b{-\f 16}k\f 1{(-80)^k}
\\&\e \Ls 5p
\sum_{k=0}^{p-1}\b{-\f 23}k\b{-\f 56}k\f 1{(-80)^k}
\\&\e \cases x-\f px\mod
{p^2}&\t{if $p\e 1,19\mod {30}$ and so $4p=x^2+75y^2$ with $3\mid
x-2$,}
\\5x-\f p{5x}\mod
{p^2}&\t{if $p\e 7,13\mod {30}$ and so $4p=25x^2+3y^2$ with $3\mid
x-1$,}
 \\0\mod p&\t{if
$p\e 17,23\mod{30}$}\endcases\endaligned$$
\endpro

\pro{Theorem 4.7} Let $p$ be an odd prime and $p\not=11$. Then

$$\aligned&\sum_{k=0}^{[p/6]}\b{[\f p3]}k\b{[\f p6]}k\Ls {27}{16}^k
\\&\e \cases -\sls{-11+x/y}p\sls x{11}x\mod
p&\t{if $\sls p3=\sls p{11}=1$ and so $4p=x^2+11y^2$,}
\\4y\sls{-11+\sls x{11}x/y}p\mod
p&\t{if $\sls p{11}=-\sls p3=1$ and so $4p=x^2+11y^2$,}
 \\0\mod p&\t{if
$\sls p{11}=-1$.}\endcases\endaligned$$
\endpro
\pro{Conjecture 4.7} Let $p>11$ be a prime such that $\sls p{11}=1$
and so $4p=x^2+11y^2$. Then

$$\aligned\sum_{k=0}^{p-1}\b{-\f 13}k\b{-\f 16}k\Ls {27}{16}^k
&\e\Big(-\f{11}{16}\Big)^{\f{1-\sls p3}2} \sum_{k=0}^{p-1}\b{-\f
23}k\b{-\f 56}k\Ls {27}{16}^k
\\& \e \cases -\sls{-11+x/y}p\sls x{11}(x-\f px)\mod {p^2}&\t{if
$3\mid p-1$,}
\\-\f 14\sls{-11+\sls x{11}x/y}p(11y-\f py)\mod
{p^2}&\t{if $3\mid p-2$}
 \endcases\endaligned$$
and so
$$\align &y\Ls{-11+\sls x{11}x/y}p\\&\e \f 14\sum_{k=0}^{(p-5)/6}
\b{\f{p-2}3}k\b{\f{p-5}6}k\Ls{27}{16}^k\mod p\qtq{for}4p=x^2+11y^2\e
2\mod 3.\endalign$$

\endpro
\pro{Theorem 4.8} Let $p>3$ be a prime. Then
$$\aligned&\sum_{k=0}^{[p/6]}\b{[\f p3]}k\b{[\f p6]}k\Big(-\f 9{16}\Big)^k
\\&\e \cases L\mod p&\t{if $p\e 1\mod 3$ and $4p=L^2+27M^2$ with $3\mid
L-2$,}\\0\mod p&\t{if $p\e 2\mod 3$.}
\endcases\endaligned$$
\endpro
Proof. From [S6, Theorem 3.2] we know that
$$P_{[\f p3]}\Ls 54\e \cases
L\mod p&\t{if $p\e 1\mod 3$ and $4p=L^2+27M^2$ with $3\mid
L-2$,}\\0\mod p&\t{if $p\e 2\mod 3$.}
\endcases$$
Thus taking $m=3$ and $t=\f{25}{16}$ in Theorem 2.1 and then
applying the above we deduce the result.

\pro{Conjecture 4.8} Let $p\e 1\mod 3$ be a prime and so
$4p=L^2+27M^2$ with $3\mid L-2$. Then
$$\sum_{k=0}^{p-1}\b{-\f 13}k\b{-\f 16}k\Big(-\f 9{16}\Big)^k
\e \sum_{k=0}^{p-1}\b{-\f 23}k\b{-\f 56}k\Big(-\f 9{16}\Big)^k\e
L-\f pL\mod {p^2}.$$\endpro

\pro{Theorem 4.9} Let $p$ be an odd prime. Then
$$\aligned&\sum_{k=0}^{[p/6]}\b{[\f p3]}k\b{[\f p6]}k\Ls{27}2^k
\\&\e \cases (-1)^{[\f p8]}\sls{-2-c/d}p2c\mod p&\t{if $p=c^2+2d^2\e
1,19\mod{24}$ and $4\mid c-1$,}
\\-\f 45(-1)^{[\f p8]}\sls{2+c/d}pd
\mod p&\t{if $p=c^2+2d^2\e 11,17\mod{24}$ and $4\mid c-1$,}
\\0\mod p&\t{if $p\e 5,7\mod 8$.}
\endcases\endaligned$$
\endpro
Proof. By [S6, Theorem 4.3],
$$\aligned&P_{[\f p3]}(5/\sqrt {-2})\\&\e\cases
(-1)^{[\f{p}8]}\sls{-2-\sqrt{-2}}p2c \mod p&\t{if $p=c^2+2d^2\e
1,3\mod 8$ and $4\mid c-1$,}\\0\mod p&\t{if $p\e 5,7\mod
8$}\endcases\endaligned$$ Thus taking $m=3$ and  $t=-\f{25}2$ in
Theorem 2.1 and then applying the above we deduce the result.

\pro{Conjecture 4.9} Let $p\e 1,3\mod 8$ be a prime and so
$p=c^2+2d^2$ with $4\mid c-1$. Then
$$\aligned&\sum_{k=0}^{p-1}\b{-\f 13}k\b{-\f 16}k\Big(\f {27}2\Big)^k
\\&\e \Big(-\f{25}2\Big)^{\f{1-\sls p3}2}\sum_{k=0}^{p-1}\b{-\f
23}k\b{-\f 56}k\Ls {27}2^k
\\&\e\cases (-1)^{[\f{p}8]}\sls{-2-c/d}p(2c-\f p{2c}) \mod {p^2}&\t{if
$p\e 1,19\mod{24}$,}
\\(-1)^{[\f{p}8]}\sls{2+c/d}p(10d-\f{5p}{4d}) \mod {p^2}
&\t{if $p\e 11,17\mod{24}$.}
\endcases\endaligned$$\endpro

 \subheading {5. Congruences for $\sum_{k=0}^{[p/8]}\b{[\f p{8}]}k
 \b{[\f {3p}8]}k(1-t)^k\mod p$}
 \pro{Theorem 5.1}
 Let $p>3$ be a prime and $t\in\Bbb Z_p$ with
$t\not\e 0\mod p$. Then
$$\aligned &\sum_{k=0}^{p-1}\b{-\f 1{8}}k\b{-\f 38}k(1-t)^k
\\&\e t^{\ag{-\f 18}}\sum_{k=0}^{p-1}\b{-\f 1{8}}k\b{-\f 58}k\Big(1-\f 1t\Big)^k
\\&\e \cases P_{[\f p4]}(\sqrt t)\mod p&\t{if $p\e
1,3\mod {8}$,}
\\\sls tp\sqrt tP_{[\f p4]}(\sqrt t)\mod p&\t{if $p\e
5,7\mod {8}$.}
\endcases\endaligned$$
\endpro
Proof. Taking $a=-\f 18$ in Theorem 2.2 and applying Lemmas 2.2 and
2.3 we see that

$$\aligned &\sum_{k=0}^{p-1}\b{-\f 1{8}}k\b{-\f 38}k(1-t)^k
\\&\e t^{\ag{-\f 18}}\sum_{k=0}^{p-1}\b{-\f 1{8}}k\b{-\f 5{8}}k\Big(1-\f
1t\Big)^k \e P_{2\ag{-\f 1{8}}}(\sqrt t)
\\&=\cases P_{\f{p-1}4}(\sqrt t)\mod p&\t{if $p\e
1\mod{8}$,}
\\P_{\f{3p-1}4}(\sqrt t)
\e P_{\f{p-3}4}(\sqrt t) \mod p&\t{if $p\e 3\mod{8}$,}
\\P_{\f{5p-1}4}(\sqrt t)\e  (\sqrt t)^pP_{\f{p-1}4}(\sqrt t)
\e \sls tp\sqrt tP_{\f{p-1}4}(\sqrt t)\mod p&\t{if $p\e 5\mod {8}$,}
\\P_{\f{7p-1}4}(\sqrt t)\e
(\sqrt t)^pP_{\f{3p-1}4}(\sqrt t) \e \sls tp\sqrt
tP_{\f{p-3}4}(\sqrt t) \mod p& \t{if $p\e 7\mod
{8}$.}\endcases\endaligned$$ Thus the theorem is proved.

\pro{Theorem 5.2} Let $p$ be a prime such that $p\e \pm 1\mod 8$.
Then
$$\aligned&\sum_{k=0}^{[p/8]}\b{[\f p8]}k\b{[\f {3p}8]}k\f 1{9^k}
\\&\e \cases 2x\sls x3\cdot 2^{\f{p-1}4}\mod p&\t{if $p=x^2+6y^2\e 1\mod{24}$,}
\\3x\sls x3\cdot 2^{\f{p-3}4}\mod p&\t{if $p=x^2+6y^2\e 7\mod{24}$,}
\\0\mod p&\t{if $p\e 17,23\mod{24}$.}
\endcases\endaligned$$
\endpro
Proof. From [S7, Theorem 3.5] we know that
$$P_{[\f p4]}\Ls{2\sqrt 2}3\e
\cases (-1)^{\f{p-1}2}\sls{\sqrt 2}p\sls x32x\mod p&\t{if
$p=x^2+6y^2\e 1,7\mod {24}$,}
\\0\mod p&\t{if $p\e 17,23\mod {24}$}\endcases $$
Taking $m=4$ and $t=\f 89$ in Theorem 2.1 and then applying the
above we deduce
$$\aligned &\sum_{k=0}^{[p/8]}\b{[\f p8]}k\b{[\f {3p}8]}k\f 1{9^k}
\\&\e\cases P_{[\f p4]}(\f{2\sqrt 2}3)\e 2x\sls x3\sls{\sqrt 2}p
\e 2x\sls x3\cdot 2^{\f{p-1}4}\mod p\\\qq\qq\q\qq \t{if
$p=x^2+6y^2\e 1\mod{24}$,}
\\\f 3{2\sqrt 2}P_{[\f p4]}(\f{2\sqrt 2}3)\e
\f 3{2\sqrt 2}\sls{\sqrt 2}p\sls x3 2x\e 3x\sls x3\cdot
2^{\f{p-3}4}\mod p
\\\qq\q\qq\qq\t{if $p=x^2+6y^2\e 7\mod {24}$,}\\0\mod p\ \;\qq\t{if
$p\e 17,23\mod{24}$.}
\endcases\endaligned$$
This proves the theorem.
 \pro{Conjecture 5.1} Let $p$ be a prime
such that $p\e 5,11\mod {24}$ and so $p=2x^2+3y^2(x,y\in\Bbb Z)$.
Then
$$\aligned\sum_{k=0}^{[p/8]}\b{[\f p8]}k\b{[\f {3p}8]}k\f 1{9^k}
\e \cases 3x\f ab\mod p&\t{if $p=a^2+b^2\e 5\mod{24}$,}
\\2x\f cd\mod p&\t{if $p=c^2+2d^2\e 11\mod{24}$.}\endcases\endaligned
$$\endpro
\pro{Lemma 5.1 ([S4, Lemma 4.1])} Let $p$ be an odd prime. Then
$$P_{[\f p4]}(\sqrt t)\e -\sum_{n=0}^{p-1}(n^3+4n^2+2(1-\sqrt
t)n)^{\f{p-1}2}\mod p.$$
\endpro
\pro{Theorem 5.3} Let $p\e 1,9\mod{20}$ be a prime and so
$p=x^2+5y^2$ with $x,y\in\Bbb Z$. Then
$$\aligned \sum_{k=0}^{[p/8]}\b{[\f p8]}k\b{[\f {3p}8]}k\f 1{(-4)^k}
\e\cases 2x\mod p&\t{if $p\e 1,9\mod{40}$,}
\\ 4ya/b\mod p&\t{if $p=a^2+b^2\e 21,29\mod{40}$.}
\endcases\endaligned$$
\endpro
Proof. From [LM, Theorem 11] we know that
$$\sum_{n=0}^{p-1}\Ls{n^3+4n^2+(2-\sqrt 5)n}p=2x.$$
Thus, taking $m=4$ and $t=\f 54$ in Theorem 2.1 and Lemma 5.1 we
obtain
$$\aligned &\sum_{k=0}^{[p/8]}\b{[\f p8]}k\b{[\f {3p}8]}k\f 1{(-4)^k}
\\&\e\cases P_{[\f p4]}(\sqrt{\f 54})\e 2x\mod p&\t{if $p\e 1,9\mod{40}$,}
\\ \f 2{\sqrt 5}P_{[\f p4]}(\sqrt{\f 54})\e \f 45\sqrt 5x\e 4y\f ab\mod p
&\t{if $p=a^2+b^2\e 21,29\mod{40}$.}
\endcases\endaligned$$
This is the result.

\pro{Theorem 5.4} Let $p$ be a prime such that $p\e
1,9,11,19\mod{40}$ and so $p=x^2+10y^2$ with $x,y\in\Bbb Z$. Then
$$\aligned 2x&\e \sum_{k=0}^{[p/8]}\b{[\f p8]}k\b{[\f {3p}8]}k\f
1{81^k}\mod p.
\endaligned$$
\endpro
Proof. From [LM] and Deuring's theorem we deduce that (see [S6])
$$\sum_{n=0}^{p-1}\Ls{n^3+4n^2+(2-\f 89\sqrt 5)n}p
= 2x.$$ Thus, taking $m=4$ and $t=\f{80}{81}$ in Theorem 2.1 and
then applying Lemma 5.1 and the above we deduce the result.

\pro{Theorem 5.5} Suppose that $p$ is a prime such that
$\sls{-1}p=\sls {13}p=1$ and so $p=x^2+13y^2(x,y\in\Bbb Z)$. Then
$$\aligned \sum_{k=0}^{[p/8]}\b{[\f p8]}k\b{[\f {3p}8]}k\f 1{(-324)^k}
\e\cases 2x\mod p&\t{if $p\e 1\mod 8$,}
\\\f{36}5y\f ab\mod p&\t{if $p=a^2+b^2\e 5\mod
8$.}\endcases\endaligned$$
\endpro
Proof. From [LM] and Deuring's theorem  we deduce (see [S6])
$$\sum_{n=0}^{p-1}\Ls{n^3+4n^2+2(1-\f 5{18}\sqrt{13})n}p
= 2x.$$
 Now taking $m=4$ and $t=\f{325}{324}$ in Theorem 2.1 and
then applying Lemma 5.1 and the above we deduce
$$\aligned&\sum_{k=0}^{[p/8]}\b{[\f p8]}k\b{[\f {3p}8]}k\f 1{(-324)^k}
\\&\e\cases P_{\f{p-1}4}\sls{5\sqrt{13}}{18}\e 2x\mod p\\\qq\qq\t{if $p\e
1\mod 8$,}
\\\f{18}{5\sqrt{13}}P_{\f{p-1}4}\sls{5\sqrt{13}}{18}
\e \f{18}{5\sqrt{13}}\cdot 2x\e \f {36x}{5\sqrt{-13}\sqrt {-1}}\e
\f{36}5y\f ab\mod p\\\qq\qq\t{if $p=a^2+b^2\e 5\mod
8$.}\endcases\endaligned$$ This completes the proof.
\newline{\bf Remark 5.1} Let $d\in\{5,10,13\}$, and
$f(d)=-4,81,-324$ according as $d=5,10,13$. Let $p$ be a prime such
that $p=x^2+dy^2\e 1\mod 8$. After reading the author's conjectures
on $$\align\sum_{k=0}^{p-1}\f{\b{-\f 13}k\b{-\f 16}k}{m^k},\q
\sum_{k=0}^{p-1}\f{\b{-\f 14}k^2}{m^k}, \q\sum_{k=0}^{p-1}\f{\b{-\f
13}k^2}{m^k},\q
 \sum_{k=0}^{p-1}\f{\b{-\f
12}k\b{-\f 14}k}{m^k}\mod {p^2},\endalign$$ the author's brother
Z.W. Sun conjectured
$$ \sum_{k=0}^{p-1}\b{-\f 1{8}}k\b{-\f 38}k\f 1{f(d)^k}
\e 2x-\f p{2x}\mod{p^2}.$$

\pro{Theorem 5.6} Suppose that $p$ is a prime such that
$\sls{-1}p=\sls {37}p=1$ and so $p=x^2+37y^2(x,y\in\Bbb Z)$. Then
$$\aligned \sum_{k=0}^{[p/8]}\b{[\f p8]}k\b{[\f {3p}8]}k\f 1{(-882^2)^k}
\e\cases 2x\mod p&\t{if $p\e 1\mod 8$,}
\\\f{1764a}{145b}y\mod p&\t{if $p=a^2+b^2\e 5\mod
8$.}\endcases\endaligned$$
\endpro
Proof. From [LM] and Deuring's theorem  we deduce (see [S6])
$$\sum_{n=0}^{p-1}\Ls{n^3+4n^2+2(1-\f {145}{882}\sqrt{37})n}p
= 2x.$$
 Now taking $m=4$ and $t=\f{37\cdot 145^2}{882^2}$ in Theorem 2.1 and
then applying Lemma 5.1 and the above we deduce
$$\aligned&\sum_{k=0}^{[p/8]}\b{[\f p8]}k\b{[\f {3p}8]}k\f 1{(-882^2)^k}
\\&\e\cases P_{\f{p-1}4}\sls{145\sqrt{37}}{882}\e
 2x\mod p\\\qq\qq\t{if $p\e
1\mod 8$,}
\\\f{882}{145\sqrt{37}}P_{\f{p-1}4}\sls{145\sqrt{37}}{882}
\e \f{882}{145\sqrt{37}}\cdot 2x\e \f {1764x}{145bx/(ay)}=
\f{1764a}{145b}y\mod p\\\qq\qq\t{if $p=a^2+b^2\e 5\mod
8$.}\endcases\endaligned$$ This completes the proof.

\pro{Conjecture 5.2} Suppose that $p$ is a prime such that
$\sls{-1}p=\sls {37}p=1$ and so $p=x^2+37y^2(x,y\in\Bbb Z)$. Then
$$ \sum_{k=0}^{p-1}\b{-\f 18}k\b{-\f 38}k\f 1{(-882^2)^k}
\e 2x-\f p{2x}\mod {p^2}.$$\endpro

\pro{Theorem 5.7} Let $p$ be a prime such that $\sls
2p=\sls{-11}p=1$ and so $p=x^2+22y^2$ with $x,y\in\Bbb Z$. Then
$$\sum_{k=0}^{[p/8]}\b{[\f p8]}k\b{[\f {3p}8]}k\f 1{9801^k}
\e\cases 2x\mod p&\t{if $p\e 1\mod 8$,}
\\\f{99}{70}\cdot 2^{\f{p+1}4}x\mod p&\t{if $p\e 7\mod 8$.}
\endcases$$
\endpro
Proof. From [LM] and Deuring's theorem we deduce (see [S6])
$$\sum_{n=0}^{p-1}\Ls{n^3+4n^2+2(1-\f {70}{99}\sqrt 2)n}p
= 2x.$$ Thus, taking $m=4$ and $t=\f{9800}{9801}$ in Theorem 2.1 and
then applying Lemma 5.1 and the above we deduce the result.

\pro{Conjecture 5.3} Let $p$ be a prime such that $\sls 2p=\sls
{-11}p=-1$ and so $p=2x^2+11y^2(x,y\in\Bbb Z)$. Then
$$\aligned\sum_{k=0}^{[p/8]}\b{[\f p{8}]}k\b{[\f {3p}8]}k\f 1{9801^k}
\e \cases 2x\f cd\mod p&\t{if $p=c^2+2d^2\e 3\mod{8}$,}
\\\f{99}{35}x\f ab\mod p&\t{if $p=a^2+b^2\e 5\mod{8}$.}\endcases\endaligned
$$\endpro

\pro{Theorem 5.8} Let $p$ be a prime such that $\sls
{-2}p=\sls{29}p=1$ and so $p=x^2+58y^2$ with $x,y\in\Bbb Z$. Then
$$ 2x \e \sum_{k=0}^{[p/8]}\b{[\f p8]}k\b{[\f {3p}8]}k\f 1{99^{4k}}
\mod p.$$
\endpro

Proof. From [LM] and Deuring's theorem  we have (see [S6])
$$\sum_{n=0}^{p-1}\Ls{n^3+4n^2+2(1-\f {1820}{99^2}\sqrt {29})n}p
= 2x.$$ Thus, taking $m=4$ and $t=\f{29\cdot 1820^2}{99^4}$ in
Theorem 2.1 and then applying Lemma 5.1 and the above we deduce the
result.

\pro{Conjecture 5.4} Let $p$ be a prime such that $\sls {-2}p=\sls
{29}p=-1$ and so $p=2x^2+29y^2(x,y\in\Bbb Z)$. Then
$$y\e \f{910}{9801}
\sum_{k=0}^{[p/8]}\b{[\f p{8}]}k\b{[\f {3p}8]}k\f 1{99^{4k}} \mod
p.$$
\endpro
\pro{Theorem 5.9} Let $p$ be a prime such that $p\e 1,19\mod{24}$
and so $p=x^2+2y^2$. Then
$$2x\e \sum_{k=0}^{[p/8]}\b{[\f p8]}k\b{[\f {3p}8]}k\f 1{2401^k}\mod p.$$
\endpro
Proof. From [S7, the proof of Theorem 3.9] we know that
$$\sum_{n=0}^{p-1}\Ls{n^3+4n^2+(2-\f {40}{49}\sqrt{6})n}p
= 2x.$$ Thus, taking $m=4$ and $t=\f{2400}{2401}$ in Theorem 2.1 and
then applying Lemma 5.1 and the above we deduce the result.

\pro{Conjecture 5.5} let $p$ be a prime such that $p\e 1,3\mod{8}$
and so $p=x^2+2y^2$. Then
$$ \sum_{k=0}^{p-1}\b{-\f 1{8}}k\b{-\f 38}k\f 1{2401^k}\e 2x-\f p{2x}
\mod {p^2}.$$
\endpro

\pro{Theorem 5.10} Let $p>7$ be a prime. Then
$$\aligned&\sum_{k=0}^{[p/8]}\b{[\f p8]}k\b{[\f {3p}8]}k\Big(\f{256}{81}\Big)^k
\\&\e\cases
2x\sls{3(7+x/y)}p\sls x7\mod p&\t{if $p=x^2+7y^2\e 1,3\mod {8}$,}
\\\f{18}{5}y\sls{3(7+\sls x7x/y)}p\mod p
&\t{if $p=x^2+7y^2\e 5,7\mod {8}$,}
\\0\mod p&\t{if $p\e 3,5,6\mod 7$.}
\endcases\endaligned$$
\endpro
Proof. From [S7, Theorem 3.1] we know that
$$P_{[\f p4]}\Ls{5\sqrt{-7}}9
\e\cases 2x\sls{3(7+x/y)}p\sls x7\mod p&\t{if $p=x^2+7y^2\e
1,2,4\mod 7$,}\\0\mod p&\t{if $p\e 3,5,6\mod 7$.}
\endcases$$
Now taking $m=4$ and $t=-\f{175}{81}$ in Theorem 2.1 and then
applying the above we deduce the result.

\pro{Theorem 5.11} Let $p>7$ be a prime. Then
$$\aligned&\sum_{k=0}^{[p/8]}\b{[\f p8]}k\b{[\f {3p}8]}k\Big(-\f{256}{3969}\Big)^k
\\&\e\cases
2\sls p3\sls x7x\mod p&\t{if $p=x^2+7y^2\e 1,3\mod {8}$,}
\\-\f{126}{65}\sls p3\sls x7x\mod p
&\t{if $p=x^2+7y^2\e 5,7\mod {8}$,}
\\0\mod p&\t{if $p\e 3,5,6\mod 7$.}
\endcases\endaligned$$
\endpro
Proof. From [S4, Theorem 2.6] we know that
$$P_{[\f p4]}\Big(\f {65}{63}\Big)\e \cases (-1)^{[\f p4]}2x\sls p3\sls x7
\mod p&\t{if $p=x^2+7y^2\e 1,2,4\mod 7$,}
\\0\mod p&\t{if $p\e 3,5,6\mod 7$.}
\endcases\tag 5.1$$
Now taking $m=4$ and $t=\f{65^2}{63^2}$ in Theorem 2.1 and then
applying (5.1) we deduce the result.

\pro{Theorem 5.12} Let $p>7$ be a prime. Then
$$\aligned&\sum_{k=0}^{[p/8]}\b{[\f p8]}k\b{[\f {3p}8]}k\Ls{32}{81}^k
\\&\e\cases \sls p32a\mod p&\t{if $p=a^2+b^2\e 1\mod 8$ and $4\mid
a-1$,}
\\\f{18}7\sls p3a\mod p&\t{if $p=a^2+b^2\e 5\mod 8$ and $4\mid
a-1$,}
\\0\mod p&\t{if $p\e 3\mod 4$.}
\endcases\endaligned$$
\endpro
Proof. From [S4, Theorem 2.4] we know that
$$P_{[\f p4]}\Big(\f 79\Big)\e \cases \sls p32a\mod
p&\t{if $p=a^2+b^2\e 1\mod 4$ and $4\mid a-1$,}
\\0\mod p&\t{if $p\e 3\mod 4$.}
\endcases\tag 5.2$$
Now taking $m=4$ and $t=\f{49}{81}$ in Theorem 2.1 and then applying
(5.2) we deduce the result.

\pro{Theorem 5.13} Let $p>7$ be a prime. Then
$$\aligned&\sum_{k=0}^{[p/8]}\b{[\f p8]}k\b{[\f {3p}8]}k\Big(-\f{16}9\Big)^k
\\&\e\cases
2A\mod p&\t{if $p=A^2+3B^2\e 1,19\mod {24}$ and $3\mid A-1$,}
\\-\f{6}5A\mod p&\t{if $p=A^2+3B^2\e 7,13\mod {24}$ and $3\mid A-1$,}
\\0\mod p&\t{if $p\e 2\mod 3$.}
\endcases\endaligned$$
\endpro
Proof. From [S4, Theorem 2.5] we know that
$$P_{[\f p4]}\Big(\f 53\Big)\e \cases (-1)^{[\f p4]}2A\mod
p&\t{if $p=A^2+3B^2\e 1\mod 3$ and $3\mid A-1$,}
\\0\mod p&\t{if $p\e 2\mod 3$.}
\endcases\tag 5.3$$
Now taking $m=4$ and $t=\f{25}{9}$ in Theorem 2.1 and then applying
(5.3) we deduce the result.

\subheading{6. Congruences for $\sum_{k=0}^{p-1}\b{-\f 13}k^2m^k$
and $\sum_{k=0}^{p-1}\b{-\f 14}k^2m^k\mod p$}

\pro{Theorem 6.1} Let $p>3$ be a prime. Then
$$\aligned&\sum_{k=0}^{p-1}\b{-\f 13}k^29^k\e \f 1{3^{[\f p3]}}
\sum_{k=0}^{p-1}\b{-\f 13}k^2\f 1{9^k}
\\&\e\cases L\mod p&\t{if $p\e 1\mod 3$, $4p=L^2+27M^2$ and
$3\mid L-2$,}
\\0\mod{p}&\t{if $p\e 2\mod 3$.}\endcases\endaligned$$\endpro
Proof. By Theorem 2.3 we have
$$\aligned &9^{\langle -\f 13\rangle_p}
\sum_{k=0}^{p-1}\b{-\f 13}k^2\f 1{9^k}\\&\e \sum_{k=0}^{p-1}\b{-\f
13}k^29^k \e 8^{\langle -\f 13\rangle_p}P_{\langle -\f
13\rangle_p}\Ls{10}8
\\&=\cases 8^{\f{p-1}3}P_{\f{p-1}3}\sls 54 \e P_{\f{p-1}3}\sls 54\mod
p &\t{if $p\e 1\mod 3$,}
\\8^{\f{2p-1}3}P_{\f{2p-1}3}\sls 54\e 2P_{\f{p-2}3}\sls 54\mod
p&\t{if $p\e 2\mod 3$.}\endcases\endaligned$$ From [S6, Theorem 3.2]
we know that
$$P_{[\f p3]}\Ls 54\e \cases L\mod p&\t{if $p\e 1\mod 3$,}
\\0\mod p&\t{if $p\e 2\mod 3$.}\endcases$$
Since
$$9^{\ag{-\f 13}}=\cases 9^{\f{p-1}3}\e 3^{-\f{p-1}3}\mod p
&\t{if $p\e 1\mod 3$,}
\\9^{\f{2p-1}3}\e 3^{-\f{p-2}3}\mod p&\t{if $p\e 2\mod 3$,}
\endcases$$ combining all the above we deduce the result.

 \pro{Conjecture 6.1} Let $p$ be a prime such that $p\e 1\mod 3$,
  $4p=L^2+27M^2(L,M\in\Bbb Z)$ and $L\e 2\mod 3$.  Then
$$\sum_{k=0}^{p-1}\b{-\f 13}k^29^k
\e L-\f pL\mod{p^2}$$\endpro

\pro{Theorem 6.2} Let $p$ be an odd prime. Then
$$\aligned&\sum_{k=0}^{p-1}\b{-\f 14}k^2(-8)^k
\\&\e \cases (-1)^{\f{p-1}4}2x\mod p&\t{if $p=x^2+y^2\e 1\mod 4$
and $4\mid x-1$,}
\\0\mod p&\t{if $p\e 3\mod 4$}\endcases\endaligned$$
and
$$\aligned&\sum_{k=0}^{p-1}\b{-\f 14}k^2\f 1{(-8)^k}
\\&\e\cases (-1)^{\f y4}2x\mod p
&\t{if $p=x^2+y^2\e 1\mod 8$ and $4\mid x-1$,}\\(-1)^{\f{y-2}4}
2y\mod p&\t{if $p=x^2+y^2\e 5\mod 8$ and $2\nmid x$,}
\\0\mod p&\t{if $p\e 3\mod 4$.}\endcases\endaligned$$
\endpro
Proof. By Lemma 2.2 and Theorem 2.3 we have
$$\aligned&\sum_{k=0}^{p-1}\b{-\f 14}k^2(-8)^k
\\&\e (-8)^{\ag{-\f 14}}\sum_{k=0}^{p-1}\b{-\f 14}k^2\f 1{(-8)^k}
\\&\e (-9)^{\ag{-\f 14}}P_{\ag{-\f 14}}\Ls 79=9^{\ag{-\f 14}}P_{\ag{-\f
14}}\Big(-\f 79\Big)
\\&=\cases 9^{\f{p-1}4}P_{\f{p-1}4}(-\f 79)\e \sls 3p
P_{\f{p-1}4}(-\f 79)\mod p&\t{if $p\e 1\mod 4$,}
\\9^{\f{3p-1}4}P_{\f{3p-1}4}(-\f 79)\e 9^{\f{3p-1}4}
P_{\f{p-3}4}(-\f 79)\mod p&\t{if $p\e 3\mod
4$.}\endcases\endaligned$$ From [S4, Theorem 2.4] we know that
$$P_{[\f p4]}\Big(-\f 79\Big)\e \cases (-1)^{\f{p-1}4}\sls p32x\mod
p&\t{if $p=x^2+y^2\e 1\mod 4$ and $4\mid x-1$,}
\\0\mod p&\t{if $p\e 3\mod 4$.}
\endcases$$
When $p\e 1\mod 4$, we have   $(-8)^{-\ag{-\f
14}}=(-8)^{-\f{p-1}4}\e (-2)^{\f{p-1}4}\mod p.$ It is well known
that (see [BEW])
$$2^{\f{p-1}4}\e \cases(-1)^{\f y4}\mod p&\t{if $p\e 1\mod 8$,}
\\(-1)^{\f{y-2}4}\f yx\mod p&\t{if $p\e 5\mod 8$.}
\endcases$$
Now combining all the above we deduce the result.

\pro{Conjecture 6.2} Let $p$ be a prime of the form $4k+1$  and so
$p=x^2+y^2$ with $4\mid x-1$. Then
$$\sum_{k=0}^{p-1}\b{-\f 14}k^2(-8)^k
\e(-1)^{\f{p-1}4}(2x-\f p{2x})\mod {p^2}$$ and
$$\aligned\sum_{k=0}^{p-1}\f{\b{-\f 14}k^2}{(-8)^k}
\e\cases (-1)^{\f y4}(2x-\f p{2x})\mod {p^2} &\t{if $p\e 1\mod
8$,}\\(-1)^{\f{y-2}4} (2y-\f p{2y})\mod {p^2}&\t{if $p\e 5\mod 8$.}
\endcases\endaligned$$\endpro

\pro{Theorem 6.3} Let $p>3$ be a prime. Then
$$\aligned&\sum_{k=0}^{p-1}\b{-\f 14}k^24^k
\e \f{3-(-1)^{\f{p-1}2}}2\Ls 2p\sum_{k=0}^{p-1}\b{-\f 14}k^2\f
1{4^k}
\\&\e \cases (-1)^{\f{p-1}4+\f{A-1}2}2A\mod p&\t{if $p=A^2+3B^2
\e 1\mod {12}$,}
\\(-1)^{\f{p+1}4}6B\mod p&\t{if $p=A^2+3B^2\e 7\mod {12}$
and $4\mid B-1$,}
\\0\mod p&\t{if $p\e 2\mod 3$}\endcases\endaligned$$
\endpro
Proof. By Lemma 2.2 and Theorem 2.3 we have
$$\aligned&\sum_{k=0}^{p-1}\b{-\f 14}k^24^k
\\&\e 4^{\ag{-\f 14}}\sum_{k=0}^{p-1}\b{-\f 14}k^2\f 1{4^k}
\\&\e 3^{\ag{-\f 14}}P_{\ag{-\f 14}}\Ls 53=(-3)^{\ag{-\f 14}}P_{\ag{-\f
14}}\Big(-\f 53\Big)
\\&=\cases (-3)^{\f{p-1}4}P_{\f{p-1}4}(-\f 53)\mod p&\t{if $p\e 1\mod 4$,}
\\(-3)^{\f{3p-1}4}P_{\f{3p-1}4}(-\f 53)\e (-3)^{-\f{p-3}4}
P_{\f{p-3}4}(-\f 53)\mod p&\t{if $p\e 3\mod
4$.}\endcases\endaligned$$ From [S4, Theorem 2.5] we know that
$$P_{[\f p4]}\Big(-\f 53\Big)\e \cases 2A\mod
p&\t{if $p=A^2+3B^2\e 1\mod 3$ and $3\mid A-1$,}
\\0\mod p&\t{if $p\e 2\mod 3$.}
\endcases$$
Hence the result is true for $p\e 2\mod 3$. \par Now assume
$p=A^2+3B^2\e 1\mod 3$ and $A\e 1\mod 3$. If $p\e 1\mod {12}$, by
[S2, p.1317] we have $3^{\f{p-1}4}\e (-1)^{\f{A-1}2}\mod p$ and
$4^{\ag{-\f 14}}=4^{\f{p-1}4}\e \ls 2p=(-1)^{\f{p-1}4}\mod p$. Hence
$$\sum_{k=0}^{p-1}\b{-\f 14}k^24^k
\e (-1)^{\f{p-1}4}\sum_{k=0}^{p-1}\b{-\f 14}k^2\f 1{4^k} \e
(-1)^{\f{p-1}4+\f{A-1}2}2A\mod p.$$ If $p\e 7\mod {12}$ and $B\e
1\mod 4$, by [S2, p.1317] we have $3^{\f{p-3}4}\e \f BA\mod p$.
Since $4^{\ag{-\f 14}}=4^{\f{3p-1}4}=2^{p-1+\f{p+1}2} \e 2\sls
2p\mod p$, by the above we get
$$\sum_{k=0}^{p-1}\b{-\f 14}k^24^k
\e 2\Ls 2p\sum_{k=0}^{p-1}\b{-\f 14}k^2\f 1{4^k} \e
(-1)^{\f{p-3}4}\f AB\cdot 2A\e (-1)^{\f{p+1}4}6B\mod p.$$ Now
combining all the above we deduce the result.

\pro{Conjecture 6.3} Let $p\e 1\mod 3$ be a prime and so
$p=A^2+3B^2$. Then
$$\aligned&\sum_{k=0}^{p-1}\b{-\f 14}k^24^k\e
\sum_{k=0}^{p-1}\b{-\f 14}k\b{-\f 12}k(-8)^k
\\&\e\cases (-1)^{\f{p-1}4+\f{A-1}2}
(2A-\f p{2A})\mod {p^2}&\t{if $p\e 1\mod {12}$,}
\\(-1)^{\f{p+1}4+\f{B-1}2}
(6B-\f p{2B})\mod {p^2}&\t{if $p\e 7\mod {12}$.}
\endcases\endaligned$$ and
$$\aligned\sum_{k=0}^{p-1}\f{\b{-\f 14}k^2}{4^k}
\e\cases (-1)^{\f{A-1}2}(2A-\f p{2A})\mod {p^2}&\t{if $p\e
1\mod{12}$,}
\\(-1)^{\f{B-1}2}(3B-\f p{4B})\mod {p^2}&\t{if $p\e
7\mod{12}$.}
\endcases\endaligned$$\endpro

\pro{Theorem 6.4} Let $p\not=2,7$ be a prime. Then
$$\aligned\sum_{k=0}^{p-1}\b{-\f 14}k^264^k
&\e\f{9-7(-1)^{\f{p-1}2}}2\Ls 2p\sum_{k=0}^{p-1}\b{-\f 14}k^2\f
1{64^k}
\\&\e\cases
(-1)^{\f{p-1}4+\f{x-1}2}2x\mod p&\t{if $p=x^2+7y^2\e 1\mod
4$,}\\(-1)^{\f{p+1}4+\f{y-1}2}42y\mod p&\t{if $p=x^2+7y^2\e 3\mod
4$,}
\\0\mod p&\t{if $p\e 3,5,6\mod 7$.}\endcases\endaligned$$\endpro

Proof. By Lemma 2.2 and Theorem 2.3 we have
$$\aligned&\sum_{k=0}^{p-1}\b{-\f 14}k^264^k
\\&\e 64^{\ag{-\f 14}}\sum_{k=0}^{p-1}\b{-\f 14}k^2\f 1{64^k}
\\&\e 63^{\ag{-\f 14}}P_{\ag{-\f 14}}\Ls {65}{63}
=(-63)^{\ag{-\f 14}}P_{\ag{-\f 14}}\Big(-\f {65}{63}\Big)
\\&=\cases (-63)^{\f{p-1}4}P_{\f{p-1}4}(-\f {65}{63})\mod p&\t{if $p\e 1\mod 4$,}
\\(-63)^{\f{3p-1}4}P_{\f{3p-1}4}(-\f {65}{63})\e (-63)^{-\f{p-3}4}
P_{\f{p-3}4}(-\f {65}{63})\mod p&\t{if $p\e 3\mod
4$.}\endcases\endaligned$$ By [S4, Theorem 2.6],
$$P_{[\f p4]}\Big(-\f{65}{63}\Big)
\e \cases 2x\sls p3\sls x7\mod p&\t{if $p=x^2+7y^2\e 1,2,4\mod 7$,}
\\0\mod p&\t{if $p\e 3,5,6\mod 7$.}
\endcases$$
Hence the result is true for $p\e 3,5,6\mod 7$.
\par Now suppose $p\e 1,2,4\mod 7$ and so $p=x^2+7y^2$ with $x,y\in
\Bbb Z$. If $p\e 1\mod 4$, by [S2, p.1317] we have $7^{\f{p-1}4} \e
(-1)^{\f{x-1}2}\sls x7\mod p$ and so
$$\align&(-63)^{\f{p-1}4}P_{\f{p-1}4}(-\f {65}{63})
\\&\e (-1)^{\f{p-1}4}\Ls 3p
\cdot (-1)^{\f{x-1}2}\Ls x7\cdot 2x\Ls p3\Ls x7
=(-1)^{\f{p-1}4+\f{x-1}2} 2x\mod p.\endalign$$ If $p\e 3\mod 4$, by
[S2, p.1317] we have $7^{\f{p-3}4} \e (-1)^{\f{y+1}2}\sls x7\f
yx\mod p$ and so
$$\align&(-63)^{-\f{p-3}4}P_{\f{p-3}4}(-\f {65}{63})
\\&\e (-1)^{\f{p-3}4}3\Ls 3p
\cdot (-1)^{\f{y+1}2}\Ls x7\f xy\cdot 2x\Ls p3\Ls x7 \e
(-1)^{\f{p+1}4+\f{y-1}2} 42y\mod p.\endalign$$ Note that
$$64^{\ag{-\f 14}}=\cases 64^{\f{p-1}4}\e \sls 2p\mod p&\t{if $p\e
1\mod 4$,}
\\64^{\f{3p-1}4}\e 8\ls 2p\mod p&\t{if $p\e 3\mod 4$.}
\endcases$$
Combining all the above we deduce the result.

\pro{Conjecture 6.4} Let $p>2$ be a prime such that $p\e 1,2,4\mod
7$ and so $p=x^2+7y^2$. Then
$$\aligned\sum_{k=0}^{p-1}\b{-\f 14}k^264^k
\e\cases \sls 2p(-1)^{\f{x-1}2}(2x-\f p{2x})\mod {p^2}&\t{if $p\e
1\mod 4$,}\\\sls 2p(-1)^{\f{y-1}2}(42y-\f {3p}{2y})\mod {p^2}&\t{if
$p\e 3\mod 4$}
\endcases\endaligned$$
and
$$\aligned\sum_{k=0}^{p-1}\f{\b{-\f 14}k^2}{64^k}
\e\cases (-1)^{\f{x-1}2}(2x-\f p{2x})\mod {p^2}&\t{if $p\e 1\mod
4$,}\\\f{3}4(-1)^{\f{y-1}2}(7y-\f p{4y})\mod {p^2}&\t{if $p\e 3\mod
4$.}
\endcases\endaligned$$\endpro

\pro{Conjecture 6.5} Let $p$ be a prime such that $p\e 5,7\mod 8$.
Then
$$\sum_{k=0}^{p-1}(-1)^k\b{-\f 14}k^2\e
\cases (-1)^{\f{x+1}2}(2x-\f p{2x})\mod {p^2} &\t{if $p=x^2+2y^2\e
1\mod 8$,}\\(-1)^{\f{y-1}2}(4y-\f p{2y})\mod {p^2} &\t{if
$p=x^2+2y^2\e 3\mod 8$},\\0\mod p&\t{if $p\e 5,7\mod 8$.}\endcases$$
\endpro
\pro{Conjecture 6.6} Let $p$ be an odd prime.
\par $(\t{\rm i})$ If $p\e 1\mod 4$ and so $p=x^2+y^2$ with $2\nmid x$,
then
$$\aligned\sum_{k=0}^{p-1}\f{\b{-\f 14}k\b{-\f 12}k}{4^k}
&\e\Ls p3\sum_{k=0}^{p-1}\b{-\f 12}k\b{-\f 16}k2^k
 \\&\e\cases (-1)^{\f{p-1}4+\f{x+1}2}(2x-\f
p{2x})\mod{p^2}&\t{if $12\mid p-1$,}
\\2y-\f p{2y}\mod{p^2}&\t{if $12\mid p-5$.}\endcases\endaligned$$
\par $(\t{\rm ii})$ If $p\e 3\mod 4$, then
$$\sum_{k=0}^{p-1}\f{\b{-\f 14}k\b{-\f 12}k}{4^k}
\e 0\mod{p^2}\qtq{and}\sum_{k=0}^{p-1}\b{-\f 12}k\b{-\f 16}k2^k\e
0\mod p.$$
\endpro

\pro{Conjecture 6.7} Let $p$ be an odd prime. Then
$$\aligned&\sum_{k=0}^{p-1}\f{\b{-\f 14}k\b{-\f 12}k}{(-3)^k}
\e (-1)^{\f{p-1}4}\sum_{k=0}^{p-1}\f{\b{-\f 14}k\b{-\f 12}k}{81^k}
\\&\e\cases 2x-\f p{2x}\mod {p^2}&\t{if $p=x^2+y^2\e 1\mod
4$ and $2\nmid x$,}
\\0\mod p&\t{if $p\e 3\mod 4$.}\endcases\endaligned$$\endpro

\pro{Conjecture 6.8} Let $p$ be an odd prime. Then
$$\aligned&\sum_{k=0}^{p-1}\f{\b{-\f 14}k\b{-\f 12}k}{(-80)^k}
\\&\e\cases 2x-\f p{2x}\mod {p^2}&\t{if $p=x^2+y^2\e \pm 1\mod
5$ and $2\nmid x$,}
\\2y-\f p{2y}\mod {p^2}&\t{if $p=x^2+y^2\e \pm 2\mod
5$ and $2\nmid x$,}
\\0\mod p&\t{if $p\e 3\mod 4$.}\endcases\endaligned$$\endpro

\pro{Conjecture 6.9} Let $p>5$ be a prime. Then
$$\aligned&\sum_{k=0}^{p-1}\b{-\f 14}k\b{-\f 12}k2^k
\\&=\cases 2x-\f p{2x}\mod {p^2}&\t{if $p=x^2+2y^2$ with $x\e 1\mod 4$,}
\\0\mod p&\t{if $ p\e 5,7\mod 8$.}
\endcases\endaligned$$
\endpro

\pro{Conjecture 6.10} Let $p>5$ be a prime. Then
$$\aligned&\sum_{k=0}^{p-1}\b{-\f 12}k\b{-\f 13}k(-3)^k
\e\sum_{k=0}^{p-1}\f{\b{-\f 12}k\b{-\f 13}k}{(-27)^k}
\e\sum_{k=0}^{p-1}\f{\b{-\f 12}k\b{-\f 23}k}{(-4)^k}
\\&\e \Ls
p5\sum_{k=0}^{p-1}\f{\b{-\f 12}k\b{-\f 13}k}{5^k}   \e
\Ls{-1}p\sum_{k=0}^{p-1}\b{-\f 12}k\b{-\f 13}k2^k
\\&\e\cases 2A-\f p{2A}\mod{p^2}&\t{if
$p=A^2+3B^2\e 1\mod 3$ with $3\mid A-1$,}\\0\mod p&\t{if $p\e 2\mod
3$.}
\endcases\endaligned$$\endpro

\pro{Conjecture 6.11} Let $p>5$ be a prime. Then
$$\aligned&\sum_{k=0}^{p-1}\f{\b{-\f 12}k\b{-\f 13}k}{(-4)^k}
\\&=\cases \sls p5 (2A-\f p{2A})\mod{p^2}\\\q\t{if
$p=A^2+3B^2\e 1\mod 3$ with $5\mid AB$ and $3\mid A-1$,} \\ \sls p5
(A+3B-\f p{A+3B})\mod{p^2}\\\q\t{if $p=A^2+3B^2\e 1\mod 3$ with
$A/B\e -1,-2\mod 5$ and $3\mid A-1$,}
\\0\mod
p\qq\t{if $p\e 2\mod 3$.}
\endcases\endaligned$$\endpro

\pro{Conjecture 6.12} Let $p>5$ be a prime. Then
$$
\sum_{k=0}^{p-1}\b{-\f 12}k\b{-\f 16}k\Big(-\f 3{125}\Big)^k
\e\cases 2A-\f p{2A}\mod {p^2}&\t{if $p=A^2+3B^2\e 1\mod 3$,}
\\0\mod p&\t{if $p\e 2\mod 3$.}
\endcases$$
\endpro

\enddocument

\Refs\widestnumber\key{BEW}
\ref\key AAR\by G. Andrews, R. Askey, R. Roy\book Special Functions,
\publ Encyclopedia Math. Appl., vol. 71,  Cambridge Univ. Press,
Cambridge, 1999\endref

 \ref\key B\by H. Bateman\book
Higher transcendental functions, Vol.II\publ McGraw-Hill Book Co.
Inc.\yr 1953\endref


 \ref \key BEW\by  B.C. Berndt, R.J. Evans and K.S.
Williams\book  Gauss and Jacobi Sums\publ John Wiley $\&$
Sons\publaddr New York\yr 1998\endref

\ref\key BM\by J. Brillhart and P. Morton \paper Class numbers of
quadratic fields, Hasse invariants of elliptic curves, and the
supersingular polynomial \jour J. Number Theory \vol 106\yr
2004\pages 79-111\endref

 \ref \key G\by H.W. Gould\book Combinatorial
Identities, A Standardized Set of Tables Listing 500 Binomial
Coefficient Summations\publ Morgantown, W. Va.\yr 1972\endref

\ref\key I\by N. Ishii\paper Trace of Frobenius endomorphism of an
elliptic curve with complex multiplication \jour Bull. Austral.
Math. Soc.\vol 70\yr 2004 \pages 125-142\endref

\ref\key JM\by A. Joux et F. Morain\paper Sur les sommes de
caract$\grave e$res li\'ees aux courbes elliptiques $\grave a$
multiplication complexe \jour J. Number Theory\vol 55\yr 1995\pages
108-128\endref

\ref\key LH\by D. H. Lee and S. G. Hahn \paper Gauss sums and
binomial coefficients \jour J. Number Theory\vol 92\yr 2002 \pages
257-271\endref

\ref\key LM\by F. Lepr$\acute {\t{e}}$vost and F. Morain \paper
Rev$\Hat {\t{e}}$tements de courbes elliptiques $\grave {\t{a}}$
multiplication complexe par des courbes hyperelliptiques et sommes
de caract$\grave {\t{e}}$res \jour J. Number Theory \vol 64\yr
1997\pages 165-182\endref


\ref\key M\by P. Morton\paper Explicit identities for invariants of
elliptic curves \jour J. Number Theory\vol 120\yr 2006\pages
234-271\endref

 \ref\key PV\by R. Padma and S. Venkataraman\paper Elliptic curves
with complex multiplication and a character sum\jour J. Number
Theory\vol 61\yr 1996\pages 274-282\endref

\ref\key PR\by J.C. Parnami and A.R. Rajwade\paper A new cubic
character sum\jour Acta Arith.\vol 40\yr 1982\pages 347-356\endref

\ref\key R1\by A.R. Rajwade \paper The Diophantine equation
$y^2=x(x^2+21Dx+112D^2)$ and the conjectures of Birch and
Swinnerton-Dyer \jour J. Austral. Math. Soc. Ser. A\vol 24\yr 1977
\pages 286-295\endref

\ref\key R2\by A.R. Rajwade \paper On a conjecture of Williams\jour
Bull. Soc. Math. Belg. Ser. B\vol 36\yr 1984\pages 1-4\endref

 \ref\key RPR\by D.B. Rishi, J.C. Parnami and A.R.
Rajwade \paper Evaluation of a cubic character sum using the
$\sqrt{-19}$ division points of the curve $y^2=x^3-2^3\cdot
19x+2\cdot 19^2$\jour J. Number Theory\vol 19\yr 1984\pages
184-194\endref


\ref\key S1\by Z.H. Sun\paper On the number of incongruent residues
of $x^4+ax^2+bx$ modulo $p$\jour J. Number Theory \vol 119\yr
2006\pages 210-241\endref

\ref\key S2\by Z.H. Sun\paper On the quadratic character of
quadratic units\jour J. Number Theory \vol 128\yr 2008\pages
1295-1335\endref

 \ref\key S3\by Z.H. Sun\paper
Congruences concerning Legendre polynomials\jour Proc. Amer. Math.
Soc.\vol 139\yr 2011\pages 1915-1929
\endref

\ref\key S4\by Z.H. Sun\paper Congruences concerning Legendre
polynomials II\finalinfo
 arXiv:1012.3898. http://arxiv.org/\newline abs/1012.3898\endref

\ref\key S5\by Z.H. Sun\paper Congruences concerning Legendre
polynomials III\finalinfo
 arXiv:1012.4234. http://arxiv.org/\newline abs/1012.4234\endref

\ref \key S6\by Z.H. Sun\paper Congruences involving
$\b{2k}k^2\b{3k}km^{-k}$,preprint\jour arXiv:1104.2789v3.
http://arxiv.org/\newline abs/1104.2789\endref

\ref \key S7\by Z.H. Sun\paper Congruences involving
$\b{2k}k^2\b{4k}{2k}m^{-k}$,preprint\jour arXiv:1104.3047.
http://arxiv.org/\newline abs/1104.3047\endref

\ref\key W\by K.S. Williams \paper Evaluation of character sums
connected with elliptic curves\jour Proc. Amer. Math. Soc.\vol 73\yr
1979\pages 291-299\endref

\endRefs
\bye